\documentclass{amsart}
\usepackage{array}
\usepackage{multirow}
\usepackage{graphicx}
\usepackage{epsfig}
\usepackage{amsfonts,amsmath}
\usepackage{amssymb}
\usepackage{latexsym}
\usepackage{epsf}
\usepackage{graphicx,color,graphics}
\usepackage{amsmath,amssymb}
\usepackage{dsfont}
\usepackage[latin1]{inputenc}  
\newtheorem{theorem}{Theorem}[section]

\newtheorem{lemma}[theorem]{Lemma}

\newtheorem{remark}{Remark}

\numberwithin{equation}{section}

\begin{document}

\title[\hfil CAUCHY THERMOELASTIC LAMINATED TIMOSHENKO PROBLEMS OF TYPE III]
{ON THE DECAY RATES FOR TWO CAUCHY THERMOELASTIC LAMINATED TIMOSHENKO PROBLEMS OF TYPE III WITH INTERFACIAL SLIP}

\pagenumbering{arabic}

\maketitle

\begin{center}
Aissa Guesmia\\
Institut Elie Cartan de Lorraine, UMR 7502, Universit\'e de Lorraine\\
3 Rue Augustin Fresnel, BP 45112, 57073 Metz Cedex 03, France\\
e-mail: aissa.guesmia@univ-lorraine.fr
\end{center}
\begin{abstract}
The subject of this paper is to study the decay of solutions for two systems of laminated Timoshenko beams
with interfacial slip in the whole space $\mathbb{R}$ subject to a thermal effect of type III acting only on one component. When the thermal effect is acting via the second or third component of the laminated Timoshenko beam (rotation angle displacement or dynamic of the slip), we prove that both systems are polynomially stable and obtain stability estimates in the $L^2$-norm of solutions and their higher order derivatives with respect of the space variable. The decay rates, as well as the absence and presence of the regularity-loss type property, depend on the regularity of the initial data and the speeds of wave propagations. However, when the thermal effect is acting via the first comoponent (transversal displacement), we introduce a new stability number $\chi$ and prove that the stability of systems is equivalent to $\chi\ne 0$. An application to a case of lower order coupling terms will be also given. To prove our results, we use the energy method in Fourier space combined with well chosen weight functions to build appropriate Lyapunov functionals.
\end{abstract}

\vskip0,1truecm
{\bf MSC2010:} 34B05, 34D05, 34H05. 

\section{Introduction}

A typical model of laminated Timoshenko beams of length $L$ and with interfacial slip based on the Timoshenko theory can be formulated by the system (see \cite{hans, hasp, lota1} for more details)
\begin{align}\label{s1+}
\begin{cases}
\rho_{1}\,\varphi_{tt} + k\,(u - \varphi_{x})_{x} +  F_1 = 0, \\
\rho_{2}\,(3v - u)_{tt} - b\,(3v - u)_{xx} - k\,(u - \varphi_{x}) + F_2 =0, \\
{\tilde \rho}_{3}\,v_{tt} - {\tilde k}_0\,v_{xx} + 3\,k\,(u - \varphi_{x}) + 4\,{\tilde\beta}\,v_t+  {\tilde F}_3 = 0,
\end{cases}
\end{align}
where the subscripts $x$ and $t$ denote the derivative with respect to space and time variables $x$ and $t$, respectively, 
$x\in ]0,L[$ and $t>0$, combininig with some initial data and boundary conditions at $x=0$ and $x=L$. All the coefficients are positive constants and denote some physical properties of beams. The terms $F_1 =F_1 (x,t)$, $F_2 =F_2 (x,t)$ and ${\tilde F}_3 ={\tilde F}_3 (x,t)$ are external forces and play the role of controls. The functions $\varphi =\varphi (x,t)$ and $u =u (x,t)$ represent, respectively, the transverse and rotation angle dispalcements, and the function $v=v(x,t)$ is proportional to the amount of slip along the interface, so the third equation in \eqref{s1+} describes the dynamics of the slip. 
\vskip0,1truecm
Using the change of variables 
\begin{align*}
\begin{cases}
\rho_3 =\frac{1}{9} {\tilde\rho}_3 ,\quad k_1 =k,\quad k_2 =b,\quad k_3 =\frac{1}{9} {\tilde k}_0 ,\quad\beta =\frac{4}{9} {\tilde\beta}, \\
w=-3v,\quad \psi = 3v-u,\quad F_3 =\frac{1}{9} {\tilde F}_3 , 
\end{cases}
\end{align*}
the system \eqref{s1+} can be rewritten as 
\begin{align}\label{s1}
\begin{cases}
\rho_1 \varphi_{tt} -k_1\,(\varphi_{x} +\psi +w)_{x} +F_1  = 0, \\
\rho_2 \psi_{tt} - k_2\,\psi_{xx} + k_1 \,(\varphi_{x} +\psi +w) + F_2 =0, \\
\rho_3 w_{tt} - k_3\,w_{xx} + k_1 (\varphi_{x} +\psi +w) +\beta w_t +F_3 = 0.
\end{cases}
\end{align}
The system \eqref{s1} is mathematically a particular case of the following more general one of Bresse-type:
\begin{align}\label{Bs1}
\begin{cases}
\rho_1\varphi_{tt} -k_1 \,(\varphi_x +\psi +lw)_{x} -{\tilde l}k_3 (w_x -{\tilde l}\varphi) +F_1= 0,\\
\rho_2 \psi_{tt} - k_2\,\psi_{xx} +k_1\,(\varphi_x +\psi +lw)  +F_2 = 0,\\
\rho_3 w_{tt} - k_3\,(w_x -{\tilde l}\varphi)_{x} + l\,k_1\, (\varphi_x +\psi +lw) +\beta w_t +F_3 = 0,
\end{cases}
\end{align}
where $l$ and ${\tilde l}$ are positive constants. The system \eqref{Bs1} coincides with \eqref{s1} when $l=1$ and ${\tilde l}=0$. When $w=F_3 =l={\tilde l}=0$, the system \eqref{Bs1} is reduced to the following Timoshenko-type system: 
\begin{align}\label{Ts1}
\begin{cases}
\rho_1 \varphi_{tt} -k_1 \,(\varphi_x +\psi)_{x}  +F_1 = 0,\\
\rho_2 \psi_{tt} - k_2\,\psi_{xx} +k_1\,(\varphi_x +\psi) +F_2 = 0.
\end{cases}
\end{align} 
\vskip0,1truecm
The systems \eqref{s1}, \eqref{Bs1} and \eqref{Ts1} were the subject of various studies in the literature during the last thirty years, tackling well-posedness and stability questions by considering different types of controls $F_j$ (dampings, memories, heat conduction effects, ...). Let us mention here some of these studies related to our objectives in this paper.
\vskip0,1truecm
For the well-posedness and stability questions in the case of bounded domains, we refere the readers to the non exhaustive list of references \cite{agg, apal1, apal2, beim, clx, cdflr, chen, feng, gues2, gues4, gues8, gms, li, lizh, lota1, lota2, lota3, must, rapo, rvma, tata, wxy}.  
\vskip0,1truecm
We notice here that \eqref{s1} was generally considered in the literature under the following restrictions: \eqref{s1} is already damped via the control $\beta w_t$ and the speeds of the wave propagations of the last two equations in \eqref{s1} are asummed to be equal; that is 
\begin{equation}\label{betaspeeds}
\beta>0 \quad\hbox{and}\quad \frac{k_2}{\rho_2}=\frac{k_3}{\rho_3} .
\end{equation} 
\vskip0,1truecm
In case of unbounded domains, the stability of \eqref{Bs1} and \eqref{Ts1} has been also treated in the literature for the last few years. In this direction, we mention the papers \cite{s6, s1, s4, s5, s3, s2} (see also the references therein), where some polynomial stability estimates for $L^2$-norm of solutions have been proved using frictional dampings, heat conduction effects or memory controls.
\vskip0,1truecm
In this paper, we investigate the decay properties of two laminated Timoshenko beam with interfacial slip in the whole space $\mathbb{R}$ and in more general forms that the restrictions \eqref{betaspeeds} are not imposed here. In addition, only one external force $F_j$ is considered and it is generated by a thermal effect of type III. Without loss of generality, the coefficients $\rho_j$ in \eqref{s1} are taken equal to $1$. The first system we consider is the following:
\begin{align}\label{s21}
\begin{cases}
\varphi_{tt} -k_1\,(\varphi_{x} +\psi +w)_{x} +\tau_1 \gamma q_{xt}  = 0, \\
\psi_{tt} - k_2\,\psi_{xx} + k_1 \,(\varphi_{x} +\psi +w) + \tau_2 \gamma q_{xt} = 0, \\
w_{tt} - k_3\,w_{xx} + k_1 (\varphi_{x} +\psi +w) +\tau_3 \gamma q_{xt}  = 0,\\
q_{tt} -k_4 q_{xx} -k_5 q_{xxt} +\gamma (\tau_1 \varphi_{xt} +\tau_2 \psi_{xt} +\tau_3 w_{xt})=0,
\end{cases}
\end{align}
where $x\in \mathbb{R}$, $t>0$, $k_j >0$, $\gamma\in \mathbb{R}^*$, $q=q(x,t)$ denotes the temperature and
\begin{equation}\label{tau123}
(\tau_1 ,\tau_2 ,\tau_3 )\in \{(1,0,0), (0,1,0), (0,0,1)\}.
\end{equation} 
The thermel dissipation in \eqref{s21} is generated by the term $-k_5 q_{xxt}$ (see \eqref{Ep21} in Section 2). In the second system of interest, the thermel dissipation is generated by the term of lower order $k_5 q_{t}$; more precisely, we consider the system
\begin{align}\label{s22}
\begin{cases}
\varphi_{tt} -k_1\,(\varphi_{x} +\psi +w)_{x} +\tau_1 \gamma q_{xt}  = 0, \\
\psi_{tt} - k_2\,\psi_{xx} + k_1 \,(\varphi_{x} +\psi +w) + \tau_2 \gamma q_{xt} = 0, \\
w_{tt} - k_3\,w_{xx} + k_1 (\varphi_{x} +\psi +w) +\tau_3 \gamma q_{xt} = 0,\\
q_{tt} -k_4 q_{xx} +k_5 q_{t} +\gamma (\tau_1 \varphi_{xt} +\tau_2 \psi_{xt} +\tau_3 w_{xt})=0.
\end{cases}
\end{align}
The systems \eqref{s21} and \eqref{s22} are subject to the initial conditions
\begin{align}\label{i21}
\begin{cases}
(\varphi ,\psi , w,q )(x,0)=(\varphi_0 ,\psi_0 , w_0 ,q_0 )(x), \\
(\varphi_t ,\psi_t , w_t , q_t )(x,0)=(\varphi_1 ,\psi_1 , w_1 ,q_1 )(x).
\end{cases}
\end{align}
\vskip0,1truecm 
The main objective of this paper is to study the stability of \eqref{s21} and \eqref{s22} and get some polynomial estimates in the $L^2$-norm of solutions and their higher order derivatives with respect to $x$. We will show that, when 
$(\tau_1 ,\tau_2 ,\tau_3 )=(1,0,0)$, both \eqref{s21} and \eqref{s22} are stable if and only if $\chi\ne 0$, where 
\begin{equation}\label{chi}
\chi:= k_3 -k_2 .
\end{equation}
However, when 
\begin{equation}\label{tau12}
(\tau_1 ,\tau_2 ,\tau_3 )\in \{(0,1,0), (0,0,1)\},
\end{equation}
the systems \eqref{s21} and \eqref{s22} are always stable with a better decay rate in case 
\begin{equation}\label{k123}
k_1 =k_2 = k_3
\end{equation}
than in the opposite one. Moreover, in case \eqref{s21}, \eqref{k123} allows to avoid the regularity restriction on the initial data known as the regularity-loss property (see \cite{4, 8, s4, 6, 10, 7}). At the end of this paper, we give an application to the case where the coupling terms between the laminated Timoshenko system and the equation of heat conduction in \eqref{s21} and \eqref{s22} 
\begin{equation}\label{hoct}
\tau_j \gamma q_{xt} \quad\hbox{and}\quad \gamma (\tau_1 \varphi_{xt} +\tau_2 \psi_{xt} +\tau_3 w_{xt})
\end{equation}
are, respectively, replaced by the following lower order ones:
\begin{equation}\label{loct}
\tau_j \gamma q_{t} \quad\hbox{and}\quad -\gamma (\tau_1 \varphi_{t} +\tau_2 \psi_{t} +\tau_3 w_{t}).
\end{equation}  
\vskip0,1truecm
Our stability results show that the effect of the heat conduction is better propagated to the whole system from the second or third equation of the laminated Timoshenko system than from the first one. The proof is based on the energy method combined with the Fourier analysis (by using the transformation in the Fourier space) and well chosen weight functions.
\vskip0,1truecm
The paper is organized as follows: in Section 2, we formulate \eqref{s21} and \eqref{s22} as a first order Cauchy system and give some preliminaries. Section 3 will be devoted to the proof of some differential identities. In Section 4, we prove our stability results. We end our paper by an application to the case \eqref{loct} in Section 5.

\section{Formulation of the problems}

We start by formulating \eqref{s21} and \eqref{s22} in an abstract first order system. To do so, we introduce the new variables
\begin{align}\label{v21}
\begin{cases}
u=\varphi_t ,\quad y=\psi_t ,\quad \theta =w_t ,\quad \eta =q_t,, \\
v=\varphi_{x} +\psi +w ,\quad z=\psi_x ,\quad\phi =w_x \quad\hbox{and}\quad\sigma =q_x  .
\end{cases}
\end{align}
Then, the systems \eqref{s21} and \eqref{s22} can be presented in the form
\begin{align}\label{e10}
\begin{cases}
v_{t} - u_{x} - y - \theta = 0, \\
u_{t} - k_1 \,v_{x} +\tau_1 \gamma\,\eta_x= 0, \\
z_{t} - y_{x} = 0, \\
y_{t} - k_2\,z_{x} + k_1 v + \tau_2 \gamma\,\eta_x = 0, \\
\phi_{t} - \theta_{x} = 0,  \\
\theta_{t} - k_3\,\phi_{x} + k_1\,v  +\tau_3 \gamma \eta_x = 0, \\
\sigma_{t} - \eta_{x} = 0, \\
\eta_{t} - k_4 \sigma_{x} + (1-k_0 )k_5 \partial_x^{k_0} \eta +\gamma (\tau_1 u_x + \tau_2 y_x +\tau_3 \theta_x) = 0,
\end{cases}
\end{align}
where $k_0 =2$ in case \eqref{s21}, and $k_0 =0$ in case \eqref{s22}.
Let $U$ and its initial data $U_0$ be given by 
\begin{align*}
U = (v,\,u,\,z,\,y,\,\phi,\,\theta,\sigma,\eta )^{T}\quad\hbox{and}\quad U_0 = 
(v,\,u,\,z,\,y,\,\phi,\,\theta,\sigma,\eta )^{T} (\cdot ,0).
\end{align*}
The system \eqref{e10} and the initial conditions \eqref{i21} are reduced to 
\begin{align}\label{e11}
\begin{cases}
U_{t} (x,t)+ A_2 U_{xx} (x,t) + A_1 U_x (x,t) + A_0 U (x,t)= 0,\\
U (x, 0)=U_{0} (x),
\end{cases}  
\end{align}
where  
\begin{equation}\label{e12}
{A}_2 U_{xx} = \left(
\begin{array}{c}
0\\
\\
0 \\
\\
0   \\
\\
0  \\
\\
0\\
\\
0\\
\\
0\\
\\
-\epsilon_0 k_5 \eta_{xx}
\end{array}
\right),\,\,
{A}_1 U_{x} = \left(
\begin{array}{c}
-u_{x} \\
\\
-k_1\,v_{x} +\tau_1 \gamma\eta_x \\
\\
-\,y_{x}   \\
\\
-\,k_2 \,z_{x} +\tau_2 \gamma\eta_x   \\
\\
-\,\theta_{x} \\
\\
-\,k_3\,\phi_{x} +\tau_3 \gamma\eta_x \\
\\
-\eta_x \\
\\
-k_4\sigma_x +\gamma (\tau_1 u_x +\tau_2 y_x +\tau_3 \theta_x ) 
\end{array}
\right) \,\,\hbox{and}\,\, A_0 U= \left(
\begin{array}{c}
-y - \theta \\
\\
0  \\
\\
0   \\
\\
k_1\,v     \\
\\
0  \\
\\
k_1\,v \\
\\
0 \\
\\
(1-\epsilon_0 )k_5\eta
\end{array}
\right)
\end{equation}
and 
\begin{align}\label{epsilon0}
\epsilon_0 =\begin{cases}
1\quad \hbox{in case}\,\,\eqref{s21},\\
0\quad \hbox{in case}\,\,\eqref{s22}.
\end{cases}  
\end{align}
\vskip0,1truecm
For a function $h :\mathbb{R}\to \mathbb{C}$, $Re\, h$, $Im\, h$, $\bar{h}$ and $\widehat{h}$
denote, respectively, the real part of $h$, the imaginary part of $h$, the conjugate of $h$ and the Fourier transformation of $h$. 
Using the Fourier transformation (with respect to the space variable $x$), \eqref{e11} can be written in the Fourier space as the following first order Cauchy system:
\begin{equation}\label{g6}
\begin{aligned}
\begin{cases}
\widehat{U}_{t} (\xi,\,t)-\xi^2\,{A}_2\widehat{U} (\xi,\,t)+ i\,\xi\,{A}_1 \widehat{U} (\xi,\,t) + A_0 \widehat{U} (\xi,\,t) = 0, 
\quad &\xi\in\mathbb{R},\,\,t>0,\\
\widehat{U}(\xi,\,0) = \widehat{U}_{0}(\xi), \quad &\xi\in\mathbb{R}.
\end{cases}
\end{aligned}
\end{equation}
The solution of \eqref{g6} is given by
\begin{equation}\label{e14}
\widehat{U}(\xi,\,t) = e^{-\,(-\xi^2\,{A}_2 +i\,\xi\,{A}_1 + A_0 )\,t}\ \widehat{U}_{0}(\xi).
\end{equation}
The energy $\widehat{E}$ associated with \eqref{g6} is defined by
\begin{equation}\label{E21}
\widehat{E}(\xi,\,t) = \frac{1}{2}\left[k_1\,\vert\widehat{v}\vert^2 + \vert\widehat{u}\vert^{2} + k_2\,\vert\widehat{z}\vert^{2} + \vert\widehat{y}\vert^{2} + k_3\vert\widehat{\phi}\vert^{2} + \vert\widehat{\theta}\vert^{2} + k_4\vert\widehat{\sigma}\vert^{2} +\vert\widehat{\eta}\vert^{2} \right].
\end{equation}
System \eqref{g6} is dissipative because
\begin{equation}\label{Ep21}
\frac{d}{dt}\widehat{E}(\xi,\,t) = -k_5\xi^{2\epsilon_0} \vert\widehat{\eta}\vert^{2} . 
\end{equation}
Indeed, the first equation in \eqref{g6} is equivalent to 
\begin{align}\label{e101}
\begin{cases}
\widehat{v}_{t} - i\xi\widehat{u} - \widehat{y} - \widehat{\theta} = 0, \\
\widehat{u}_{t} - ik_1 \xi \,\widehat{v} +i\tau_1 \gamma\,\xi\widehat{\eta} = 0, \\
\widehat{z}_{t} - i\xi\widehat{y} = 0, \\
\widehat{y}_{t} - ik_2\xi\,\widehat{z} + k_1 \widehat{v} +i\tau_2 \gamma\,\xi\widehat{\eta} = 0, \\
\widehat{\phi}_{t} - i\xi\widehat{\theta} = 0,  \\
\widehat{\theta}_{t} - ik_3\,\xi\widehat{\phi} + k_1\,\widehat{v}+i\tau_3 \gamma\,\xi\widehat{\eta} = 0,\\
\widehat{\sigma}_{t} - i\xi\widehat{\eta} = 0,  \\
\widehat{\eta}_{t} -ik_4 \xi \widehat{\sigma}+k_5\xi^{2\epsilon_0}\widehat{\eta} +i\gamma\xi (\tau_1 \widehat{u} +\tau_2 \widehat{y}+\tau_3 \widehat{\theta}).
\end{cases}
\end{align}
To get \eqref{Ep21}, we have just to multiply the equations in \eqref{e101} by $k_1\bar{\widehat{v}}$, $\bar{\widehat{u}}$, $k_2\bar{\widehat{z}}$, $\bar{\widehat{y}}$, $k_3\bar{\widehat{\phi}}$, $\bar{\widehat{\theta}}$, $k_4\bar{\widehat{\sigma}}$ and $\bar{\widehat{\eta}}$, respectively, adding the obtained equations, taking the real part of the resulting expression and using the following classical relation, for two differentiable functions $h ,\,d :\mathbb{R}\to \mathbb{C}$: 
\begin{equation}\label{hd1}
\frac{d}{dt} Re\,(h {\bar d}) = Re\,(h_t {\bar d} +d_t {\bar h}).
\end{equation}
\vskip0,1truecm
We observe that the energy $\widehat{E}$ is equivalent to $\vert\widehat{U}\vert^2$ defined by
\begin{equation*}
\vert\widehat{U}(\xi,\,t)\vert^2 = \vert\widehat{v}\vert^2 + \vert\widehat{u}\vert^{2} + \vert\widehat{z}\vert^{2} + \vert\widehat{y}\vert^{2} + \vert\widehat{\phi}\vert^{2} + \vert\widehat{\theta}\vert^{2} + \vert\widehat{\sigma}\vert^{2} + \vert\widehat{\eta}\vert^{2} 
\end{equation*}
because, for $\alpha_1 =\frac{1}{2}\min\{k_1 ,k_2 ,k_3 ,k_4 ,1\}$ and $\alpha_2 =\frac{1}{2}\max\{k_1 ,k_2 ,k_3 ,k_4 ,1\}$, we have 
\begin{equation}\label{equiv}
\alpha_1 \vert\widehat{U}(\xi,\,t)\vert^2 \leq \widehat{E}(\xi,\,t) \leq \alpha_2 \vert\widehat{U}(\xi,\,t)\vert^2 ,\quad \forall \xi \in \mathbb{R}, \,\,\forall t\in \mathbb{R}_+ .
\end{equation} 
\vskip0,1truecm
Before presenting and proving our stability results in the next three sections, we prove these two lemmas that will be used in the proofs.
\vskip0,1truecm
\begin{lemma}
Let $r_1$, $r_2$ and $r_3$ be real numbers such that $r_1 >-1$ and $r_2 ,\,r_3 >0$. Then there exists $C_{r_1 ,r_2 ,r_3}>0$ such that
\begin{equation}\label{g37+}
\int_{0}^{1} \xi^{r_1}\,e^{-\,r_3\,t\,\xi^{r_2}}\ d\xi \leq C_{r_1,r_2,r_3}\,(1 + t)^{-\,(r_1 + 1)/r_2} ,\quad\forall \,t\in \mathbb{R}_+ .
\end{equation}
\end{lemma}
\vskip0,1truecm
\begin{proof} For $0\leq t\leq 1$, \eqref{g37+} is evident, for any 
$C_{r_1,r_2,r_3}\geq \frac{2^{(r_1 + 1)/r_2}}{r_1 +1}$. For $t>1$, we have 
\begin{equation*}
\int_{0}^{1} \xi^{r_1}\,e^{-\,r_3\,t\,\xi^{r_2}}\ d\xi = \int_{0}^{1} \xi^{r_1 +1-r_2}\,e^{-\,r_3\,t\,\xi^{r_2}}\,\xi^{r_2-1}\, d\xi = \int_{0}^{1} (\xi^{r_2})^{\left(r_1 + 1-r_2\right)/r_2}\,e^{-\,r_3\,t\,\xi^{r_2}}\,\xi^{r_2 -1}\ d\xi.
\end{equation*}
Taking $\tau = r_3\,t\,\xi^{r_2}$. Then 
\begin{equation*}
\xi^{r_2} = \frac{\tau}{r_3\,t} \quad\hbox{and}\quad \xi^{r_2 -1}\ d\xi = \frac{1}{r_2\,r_3\,t}\ d\tau.
\end{equation*}
Substituting in the above integral, we find
\begin{equation*}
\int_{0}^{1} (\xi^{r_2})^{\left(r_1 + 1-r_2\right)/r_2}\,e^{-\,r_3\,t\,\xi^{r_2}}\,\xi^{r_2 -1}\, d\xi = \int_{0}^{r_3\,t}\left (\frac{\tau}{r_3\,t}\right)^{\left(r_1 +1-r_2\right)/r_2}\,e^{-\,\tau}\,\frac{1}{r_2\,r_3\,t}\, d\tau 
\end{equation*}
\begin{equation*}
\leq \frac{1}{r_2\,(r_3\,t)^{\left(r_1 + 1\right)/r_3}}\int_{0}^{+\infty} \tau^{\left(r_1 + 1-r_2\right)/r_2}\,e^{-\,\tau}\ d\tau \leq \frac{2^{(r_1 + 1)/r_2}}{r_2\,r_3^{\left(r_1 + 1\right)/r_2}}C_{r_1,r_2}\,(t + 1)^{-\left(r_1 + 1\right)/r_2},
\end{equation*}
where 
\begin{equation*}
C_{r_1,r_2}=\int_{0}^{+\infty} \tau^{\left(r_1 + 1-r_2\right)/r_2}\,e^{-\,\tau}\ d\tau,
\end{equation*}
which is a convergent integral, for any $r_1 > -1$ and $r_2 >0$. This completes the proof of \eqref{g37+} with
\begin{equation*}
C_{r_1,r_2,r_3}=\max\left\{\frac{2^{(r_1 + 1)/r_2}}{r_1 +1}, \frac{2^{(r_1 + 1)/r_2}}{r_2\,r_3^{\left(r_1 + 1\right)/r_2}}C_{r_1,r_2}\right\} .
\end{equation*} 
\end{proof}
\begin{lemma}
For any positive real numbers $\sigma_1,\, \sigma_{2}$ and $\sigma_{3}$, we have
\begin{equation}\label{ineq}
\sup_{|\xi|\geq 1}\vert\xi\vert^{-\sigma_1}\,e^{-\,\sigma_{2}\,t\,\vert\xi\vert^{-\sigma_3}} \leq \left(1+\sigma_{1}/(\sigma_{2}\sigma_{3})\right)^{\sigma_{1}/\sigma_{3}}
\left(1 + t\right)^{-\,\sigma_{1}/\sigma_3},\quad\forall t\in \mathbb{R}_+ .
\end{equation}
\end{lemma}
\vskip0,1truecm
\begin{proof} 
It is clear that \eqref{ineq} is satisfied for $t=0$. Let $t>0$ and $h(x)=x^{-\sigma_1}\,e^{-\,\sigma_{2}\,t\,x^{-\sigma_3}}$, for $x\geq 1$. Direct and simple computations show that
\begin{equation*}
h^{\prime}(x)=(\sigma_2 \sigma_3 t x^{-\sigma_3} -\sigma_{1})x^{-\sigma_1 -1}\,e^{-\,\sigma_{2}\,t\,x^{-\sigma_3}}.
\end{equation*} 
If $t\geq \sigma_{1}/(\sigma_2 \sigma_3 )$, then 
\begin{eqnarray*}
h(x) &\leq & h(((\sigma_2 \sigma_3 t)/\sigma_1)^{1/\sigma_3 })=((\sigma_2 \sigma_3 )/\sigma_1 )^{-\sigma_1/\sigma_3 } e^{-\sigma_1/\sigma_3} \left(1 + 1/t\right)^{\sigma_{1}/\sigma_3}\left(1 + t\right)^{-\,\sigma_{1}/\sigma_3} \nonumber \\
& \leq & ((\sigma_2 \sigma_3 )/\sigma_1)^{-\sigma_1/\sigma_3 }\left(1 + (\sigma_2 \sigma_3 )/\sigma_1\right)^{\sigma_{1}/\sigma_3} \left(1 + t\right)^{-\,\sigma_{1}/\sigma_3} =\left(1+\sigma_{1}/(\sigma_{2}\sigma_{3})\right)^{\sigma_{1}/\sigma_{3}}
\left(1 + t\right)^{-\,\sigma_{1}/\sigma_3}, 
\end{eqnarray*}
which gives \eqref{ineq} by taking $x=|\xi|$. If $0<t< \sigma_1/(\sigma_2 \sigma_3 )$, then  
\begin{equation*}
h(x)\leq h(1)=e^{-\sigma_2 t} \left(1 + t\right)^{\sigma_{1}/\sigma_3}\left(1 + t\right)^{-\,\sigma_{1}/\sigma_3} \leq
(1+\sigma_1/(\sigma_2 \sigma_3 ))^{\sigma_1/\sigma_3 } \left(1 + t\right)^{-\,\sigma_{1}/\sigma_3}, 
\end{equation*}
which implies \eqref{ineq}, for $x=|\xi|$. 
\end{proof}

\section{Preliminary differential identities}

This section is dedicated to the proof of several identities, which will play a crucial role in the proofs. In the rest of the paper, $C$ and ${\tilde C}$ denote generic positive constants, and $C_{\varepsilon}$ denotes a generic positive constant depending on some positive constant $\varepsilon$. These generic constants can be different from line to line.  
\vskip0,1truecm
Multiplying \eqref{e101}$_{4}$ and \eqref{e101}$_{3}$ by $i\,\xi\,\overline{\widehat{z}}$ and $-i\,\xi\,\overline{\widehat{y}}$, 
respectively, adding the resulting equations, taking the real part and using \eqref{hd1}, we obtain
\begin{equation}\label{eq31}
\frac{d}{dt}Re\left(i\,\xi\,\widehat{y}\, \overline{\widehat{z}}\right) = \xi^{2}\left(\vert\widehat{y}\vert^{2} -k_2 \vert\widehat{z}\vert^{2}\right) - k_1\,Re\left(i\,\xi\,\widehat{v}\,\overline{\widehat{z}}\right) +\tau_2\gamma \xi^2\,Re\left(\widehat{\eta}\, \overline{\widehat{z}}\right).
\end{equation}
Multiplying \eqref{e101}$_{2}$ and \eqref{e101}$_{1}$ by $i\,\xi\,\overline{\widehat{v}}$ and $-i\,\xi\,\overline{\widehat{u}}$, 
respectively, adding the resulting equations, taking the real part and using \eqref{hd1}, we find
\begin{equation}\label{eq32}
\frac{d}{dt}Re\left(i\,\xi\,\widehat{u}\, \overline{\widehat{v}}\right) = \xi^{2}\left(\vert\widehat{u}\vert^{2} -k_1 \vert\widehat{v}\vert^{2}\right) - Re\left(i\,\xi\,\widehat{y}\,\overline{\widehat{u}}\right) 
- Re\left(i\,\xi\,\widehat{\theta}\,\overline{\widehat{u}}\right)+\tau_1\gamma\xi^2 \,Re\left(\widehat{\eta}\, \overline{\widehat{v}}\right).
\end{equation}
Multiplying \eqref{e101}$_{6}$ and \eqref{e101}$_{5}$ by $i\,\xi\,\overline{\widehat{\phi}}$ and $-i\,\xi\,\overline{\widehat{\theta}}$, respectively, adding the resulting equations, taking the real part and using \eqref{hd1}, we get
\begin{equation}\label{eq33}
\frac{d}{dt}Re\left(i\,\xi\,\widehat{\theta}\, \overline{\widehat{\phi}}\right) = \xi^{2}\left(\vert\widehat{\theta}\vert^{2} -k_3 \vert\widehat{\phi}\vert^{2}\right) -k_1\, Re\left(i\,\xi\,\widehat{v}\,\overline{\widehat{\phi}}\right) 
+\tau_3\gamma\xi^2 \,Re\left(\widehat{\eta}\, \overline{\widehat{\phi}}\right).
\end{equation}
Multiplying \eqref{e101}$_{6}$ and \eqref{e101}$_{1}$ by $-\xi^2\,\overline{\widehat{v}}$ and $-\xi^2\,\overline{\widehat{\theta}}$, 
respectively, adding the resulting equations, taking the real part and using \eqref{hd1}, we have
\begin{eqnarray}\label{eq34}
\frac{d}{dt}Re\left(-\xi^2\,\widehat{\theta}\, \overline{\widehat{v}}\right) & = & \xi^{2}\left(k_1\,\vert\widehat{v}\vert^{2} - \vert\widehat{\theta}\vert^{2}\right) - \xi^2\,Re\left(i\,\xi\,\widehat{u}\,\overline{\widehat{\theta}}\right) 
- k_3\,\xi^2\,Re\left(i\,\xi\,\widehat{\phi}\,\overline{\widehat{v}}\right)\nonumber \\
& & -\xi^2\,Re\left(\widehat{y}\, \overline{\widehat{\theta}}\right) +\tau_3\gamma\,\xi^2 \,Re\left(i\xi\widehat{\eta}\, \overline{\widehat{v}}\right).
\end{eqnarray}
Multiplying \eqref{e101}$_{4}$ and \eqref{e101}$_{1}$ by $-\xi^2\,\overline{\widehat{v}}$ and $-\xi^2\,\overline{\widehat{y}}$, 
respectively, adding the resulting equations, taking the real part and using \eqref{hd1}, we infer that
\begin{eqnarray}\label{eq35}
\frac{d}{dt}Re\left(-\xi^2\,\widehat{y}\, \overline{\widehat{v}}\right) & = & \xi^{2}\left(k_1\,\vert\widehat{v}\vert^{2} - \vert\widehat{y}\vert^{2}\right) - \xi^2\,Re\left(i\,\xi\,\widehat{u}\,\overline{\widehat{y}}\right) 
- k_2\,\xi^2\,Re\left(i\,\xi\,\widehat{z}\,\overline{\widehat{v}}\right) \nonumber \\
& & -\xi^2\,Re\left(\widehat{\theta}\, \overline{\widehat{y}}\right)+\tau_2\gamma\,\xi^2 \,Re\left(i\xi\widehat{\eta}\, \overline{\widehat{v}}\right).
\end{eqnarray}
Multiplying \eqref{e101}$_{8}$ and \eqref{e101}$_{7}$ by $i\,\xi\,\overline{\widehat{\sigma}}$ and $-i\,\xi\,\overline{\widehat{\eta}}$, respectively, adding the resulting equations, taking the real part and using \eqref{hd1}, we get
\begin{equation}\label{eq31p}
\frac{d}{dt}Re\left(i\,\xi\,\widehat{\eta}\, \overline{\widehat{\sigma}}\right) = \xi^{2}\left(\vert\widehat{\eta}\vert^{2} -k_4 \vert\widehat{\sigma}\vert^{2}\right) -k_5\, \xi^{2\epsilon_0} Re\left(i\,\xi\,\widehat{\eta}\,\overline{\widehat{\sigma}}\right) +\gamma\xi^2 \,Re\left(\overline{\widehat{\sigma}}\, (\tau_1 {\widehat{u}} +\tau_2 {\widehat{y}}+\tau_3 {\widehat{\theta}})\right).
\end{equation}
Multiplying \eqref{e101}$_{3}$ and \eqref{e101}$_{6}$ by $i\,\xi\,\overline{\widehat{\theta}}$ and $-i\,\xi\,\overline{\widehat{z}}$, respectively, adding the resulting equations, taking the real part and using \eqref{hd1}, we entail
\begin{equation}\label{eq36}
\frac{d}{dt}Re\left(i\,\xi\,\widehat{z}\, \overline{\widehat{\theta}}\right) = -\xi^{2}\, Re\left(\widehat{y}\,\overline{\widehat{\theta}}\right)+k_3\,\xi^{2}\, Re\left(\widehat{\phi}\,\overline{\widehat{z}}\right)
+k_1\,Re\left(i\,\xi\,\widehat{v}\,\overline{\widehat{z}}\right)-\tau_3\gamma \xi^2\,Re\left(\widehat{\eta}\, \overline{\widehat{z}}\right).
\end{equation}
Mltiplying \eqref{e101}$_{5}$ and \eqref{e101}$_{4}$ by $i\,\xi\,\overline{\widehat{y}}$ and $-i\,\xi\,\overline{\widehat{\phi}}$, respectively, adding the resulting equations, taking the real part and using \eqref{hd1}, we arrive at
\begin{equation}\label{eq37}
\frac{d}{dt}Re\left(i\,\xi\,\widehat{\phi}\, \overline{\widehat{y}}\right) = -\xi^{2}\, Re\left(\widehat{\theta}\,\overline{\widehat{y}}\right)+k_2\,\xi^{2}\, Re\left(\widehat{z}\,\overline{\widehat{\phi}}\right)
+k_1\,Re\left(i\,\xi\,\widehat{v}\,\overline{\widehat{\phi}}\right)
-\tau_2\gamma \xi^2\,Re\left(\widehat{\eta}\, \overline{\widehat{\phi}}\right).
\end{equation}
Multiplying \eqref{e101}$_{2}$ and \eqref{e101}$_{3}$ by $-\,\overline{\widehat{z}}$ and $-\,\overline{\widehat{u}}$, 
respectively, adding the resulting equations, taking the real part and using \eqref{hd1}, it follows that
\begin{equation}\label{eq310}
\frac{d}{dt}Re\left(-\,\widehat{u}\, \overline{\widehat{z}}\right) = -k_1\, Re\left(i\,\xi\,\widehat{v}\,\overline{\widehat{z}}\right) - \,Re\left(i\,\xi\,\widehat{y} \,\overline{\widehat{u}}\right)+\tau_1\gamma \,Re\left(i\xi\widehat{\eta}\, \overline{\widehat{z}}\right).
\end{equation}
Finally, multiplying \eqref{e101}$_{2}$ and \eqref{e101}$_{5}$ by $-\,\overline{\widehat{\phi}}$ and $-\,\overline{\widehat{u}}$, respectively, adding the resulting equations, taking the real part and using \eqref{hd1}, it appears that
\begin{equation}\label{eq312}
\frac{d}{dt}Re\left(-\,\widehat{u}\, \overline{\widehat{\phi}}\right) = -k_1\,Re\left(i\,\xi\,\widehat{v}\,\overline{\widehat{\phi}}\right) - \,Re\left(i\,\xi\,\widehat{\theta} \,\overline{\widehat{u}}\right)+\tau_1\gamma \,Re\left(i\xi\widehat{\eta}\, \overline{\widehat{\phi}}\right).
\end{equation} 

\section{Stability}

In this section, we investigate the asymptotic behavior, when time $t$ goes to infinity, of the solution $U$ of \eqref{e11}. First, we will show that $\vert\widehat{U}\vert^{2}$ converges exponentially to zero (with respect to time $t$) in case 
\eqref{tau12}, and in case $(\tau_1, \tau_2 ,\tau_3 )=(1,0,0)$ with $\chi \ne 0$. In case $(\tau_1, \tau_2 ,\tau_3 )=(1,0,0)$ with $\chi = 0$, we prove that $\vert\widehat{U}\vert^{2}$ does not converge to zero when $t$ goes to infinity. Let us distinguish the three cases \eqref{tau123}. 

\subsection{Case 1: $(\tau_1 ,\tau_2 ,\tau_3)=(1,0,0)$}

We start by presenting the exponential stability result for \eqref{g6} in the next lemma. 
\vskip0,1truecm
\begin{lemma}\label{lemma32}
Assume that $\chi \ne 0$; that is $k_2\ne k_3$. Let $\widehat{U}$ be a solution of \eqref{g6}. Then, there exist $c,\,\widetilde{c}>0$ such that 
\begin{equation}\label{11}
\vert\widehat{U}(\xi,\,t)\vert^{2} \leq \widetilde{c}\,e^{-\,c\,f(\xi)\,t}\,\vert\widehat{U}_0 (\xi)\vert^{2},\quad \forall \,\xi\in \mathbb{R},\quad \forall \,t\in \mathbb{R}_+ ,
\end{equation}
where  
\begin{align}\label{22}
f(\xi) =\frac{\xi^{4+2\epsilon_0 }}{{\tilde f}(\xi)}\quad\hbox{and}\quad {\tilde f}(\xi) =\begin{cases}
1 + \xi^{2} +\xi^{4} +\xi^{6} +\xi^{8}\quad &\hbox{in case}\,\,\eqref{s21},\\
1 + \xi^{2} +\xi^{4} +\xi^{6} \quad &\hbox{in case}\,\,\eqref{s22}.
\end{cases}
\end{align}
\end{lemma}
\vskip0,1truecm
\begin{proof}
Multiplying \eqref{e101}$_{2}$ and \eqref{e101}$_{8}$ by $i\frac{\vert\gamma\vert}{\gamma}\,\xi\,\overline{\widehat{\eta}}$ and $-i\frac{\vert\gamma\vert}{\gamma}\,\xi\,\overline{\widehat{u}}$, respectively, adding the resulting equations, taking the real part and using \eqref{hd1}, we get
\begin{eqnarray}\label{equ1}
\frac{d}{dt}Re\left(i\frac{\vert\gamma\vert}{\gamma}\,\xi\,\widehat{u}\, \overline{\widehat{\eta}}\right) &=& \vert\gamma\vert\xi^{2}\left(\vert\widehat{\eta}\vert^{2} - \vert\widehat{u}\vert^{2}\right) + \frac{\vert\gamma\vert}{\gamma}k_4\xi^2\,Re\left(\widehat{\sigma}\,\overline{\widehat{u}}\right) \nonumber \\ 
&& -\frac{\vert\gamma\vert}{\gamma}k_1 \xi^2\,Re\left(\widehat{v}\, \overline{\widehat{\eta}}\right) +\frac{\vert\gamma\vert}{\gamma}k_5 \xi^{2\epsilon_0}\,Re\left(i\xi\widehat{\eta}\, \overline{\widehat{u}}\right).
\end{eqnarray}
Multiplying \eqref{e101}$_{6}$ and \eqref{e101}$_{8}$ by $\overline{\widehat{\eta}}$ and $\overline{\widehat{\theta}}$, respectively, adding the resulting equations, taking the real part and using \eqref{hd1}, we find
\begin{equation}\label{equ2}
\frac{d}{dt}Re\left(\widehat{\eta}\, \overline{\widehat{\theta}}\right) = \gamma Re\left(i\,\xi\,\widehat{\theta}\,\overline{\widehat{u}}\right)+k_4\,Re\left(i\xi\widehat{\sigma}\, \overline{\widehat{\theta}}\right) -k_5 \xi^{2\epsilon_0}\,Re\left(\widehat{\eta}\, \overline{\widehat{\theta}}\right) +k_3\,Re\left(i\xi\widehat{\phi}\, \overline{\widehat{\eta}}\right)-k_1 Re\left(\widehat{v}\, \overline{\widehat{\eta}}\right).
\end{equation}
Also, multiplying \eqref{e101}$_{4}$ and \eqref{e101}$_{8}$ by $\overline{\widehat{\eta}}$ and $\overline{\widehat{y}}$, respectively, adding the resulting equations, taking the real part and using \eqref{hd1}, we obtain
\begin{equation}\label{equ3}
\frac{d}{dt}Re\left(\widehat{\eta}\, \overline{\widehat{y}}\right) = -\gamma Re\left(i\,\xi\,\widehat{u}\,\overline{\widehat{y}}\right)+k_4\,Re\left(i\xi\widehat{\sigma}\, \overline{\widehat{y}}\right) -k_5 \xi^{2\epsilon_0}\,Re\left(\widehat{\eta}\, \overline{\widehat{y}}\right) +k_2\,Re\left(i\xi\widehat{z}\, \overline{\widehat{\eta}}\right)-k_1 Re\left(\widehat{v}\, \overline{\widehat{\eta}} \right).
\end{equation}
Multiplying \eqref{e101}$_{1}$ and \eqref{e101}$_{7}$ by $i\xi\overline{\widehat{\sigma}}$ and $-i\xi\overline{\widehat{v}}$, respectively, adding the resulting equations, taking the real part and using \eqref{hd1}, we infer that
\begin{equation}\label{equ1p}
\frac{d}{dt}Re\left(i\xi\widehat{v}\, \overline{\widehat{\sigma}}\right) = -\xi^2 Re\left(\widehat{\sigma}\,\overline{\widehat{u}}\right)+\xi^2 Re\left(\widehat{v}\, \overline{\widehat{\eta}}\right) + Re\left(i\xi\widehat{y}\, \overline{\widehat{\sigma}}\right) +Re\left(i\xi\widehat{\theta}\, \overline{\widehat{\sigma}}\right).
\end{equation}
Similarily, multiplying \eqref{e101}$_{3}$ and \eqref{e101}$_{7}$ by $-\overline{\widehat{\sigma}}$ and $-\overline{\widehat{z}}$, respectively, adding the resulting equations, taking the real part and using \eqref{hd1}, we arrive at
\begin{equation}\label{equ2p}
\frac{d}{dt}Re\left(-\widehat{\sigma}\, \overline{\widehat{z}}\right) = Re\left(i\xi \widehat{\sigma}\,\overline{\widehat{y}}\right)+Re\left(i\xi\widehat{z}\, \overline{\widehat{\eta}}\right).
\end{equation}
Multiplying \eqref{e101}$_{5}$ and \eqref{e101}$_{7}$ by $-\overline{\widehat{\sigma}}$ and $-\overline{\widehat{\phi}}$, respectively, adding the resulting equations, taking the real part and using \eqref{hd1}, we entail
\begin{equation}\label{equ3p}
\frac{d}{dt}Re\left(-\widehat{\sigma}\, \overline{\widehat{\phi}}\right) = Re\left(i\xi \widehat{\sigma}\,\overline{\widehat{\theta}}\right)+Re\left(i\xi\widehat{\phi}\, \overline{\widehat{\eta}}\right).
\end{equation}
Let $\lambda_{0} ,\,\cdots,\,\lambda_{5}$ be positive constants to be defined later, and let (observe that $\chi\ne 0$ by assumption) 
\begin{equation*}
\lambda_{6} =\frac{k_2}{\chi}(\lambda_{4} +\lambda_{5}),\quad \lambda_{7} =-\frac{k_3}{\chi}(\lambda_{4} +\lambda_{5}), \quad \lambda_{8} =\frac{k_2}{k_1}\lambda_{5} \xi^2 -\lambda_{1} +\frac{k_2}{\chi}(\lambda_{4} +\lambda_{5})
\end{equation*}
and
\begin{equation*}
\lambda_{9} =\frac{k_3}{k_1}\lambda_{4} \xi^2 -\lambda_{3} -\frac{k_3}{\chi}(\lambda_{4} +\lambda_{5}).
\end{equation*}
We define the functional ${F}_0$ as follows:
\begin{eqnarray}\label{F0}
{F}_0 (\xi,\,t) & = & Re\left[i\,\xi\,\left(\lambda_1\,\widehat{y}\, \overline{\widehat{z}}+ \lambda_2\,\widehat{u}\, \overline{\widehat{v}} +\lambda_3\widehat{\theta}\, \overline{\widehat{\phi}} +\widehat{\eta}\, \overline{\widehat{\sigma}} +\lambda_6\widehat{z}\, \overline{\widehat{\theta}}+\lambda_7\widehat{\phi}\,\overline{\widehat{y}}\right)\right] \nonumber \\ 
& & +Re\left(-\lambda_4\xi^2 \,\widehat{\theta}\, \overline{\widehat{v}} -\lambda_5\xi^2 \,\widehat{y}\, \overline{\widehat{v}}-\lambda_8 \widehat{u}\, \overline{\widehat{z}}
-\lambda_{9} \widehat{u}\, \overline{\widehat{\phi}}\right).
\end{eqnarray}
By multiplying \eqref{eq31}-\eqref{eq312} by $\lambda_1 ,\,\cdots ,\,\lambda_5$, $1$ and 
$\lambda_6 ,\,\cdots,\,\lambda_{9}$, respectively, and adding the obtained equations, we see that,
thanks to the choices of $\lambda_6 ,\,\cdots ,\,\lambda_{9}$, the expression of $\frac{d}{dt}{F}_0$ does not contain the terms $Re\left(i\,\xi\,\widehat{v}\, \overline{\widehat{z}}\right)$,
$Re\left(i\,\xi\,\widehat{v}\, \overline{\widehat{\phi}}\right)$, 
$Re\left(\widehat{y}\, \overline{\widehat{\theta}}\right)$ and $Re\left(\widehat{\phi}\, \overline{\widehat{z}}\right)$ because their coefficients vanish. So, we find
\begin{eqnarray}\label{g26+}
\frac{d}{dt}{F}_0 (\xi,\,t) &=& - \xi^{2}\left(k_3\lambda_3 \,\vert\widehat{\phi}\vert^{2} +\left(\lambda_5 -\lambda_1 \right)\,\vert\widehat{y}\vert^{2} +\left(\lambda_4 -\lambda_3 \right)\,\vert\widehat{\theta}\vert^{2} +(k_1\lambda_2 -k_1\lambda_4 -k_1 \lambda_5 )\,\vert\widehat{v}\vert^{2}\right)  \nonumber \\ 
&& - \xi^{2}\left(k_2\lambda_1\,\vert\widehat{z}\vert^{2} 
+k_4\,\vert\widehat{\sigma}\vert^{2} \right)
+I_1 Re (i\xi\widehat{\theta}\, \overline{\widehat{u}})+I_2 Re (i\xi\widehat{y}\, \overline{\widehat{u}}) +
\gamma\xi^{2} Re (\widehat{\sigma}\, \overline{\widehat{u}}) \nonumber \\ 
&& +\xi^{2} (\lambda_2\vert\widehat{u}\vert^{2}
+\vert\widehat{\eta}\vert^{2}) +Re\left( i\gamma\lambda_8\xi\,\widehat{\eta}\, \overline{\widehat{z}} + i\gamma\lambda_{9}\xi\widehat{\eta}\, \overline{\widehat{\phi}} -ik_5 \xi^{2\epsilon_0 +1}\,\widehat{\eta}\, \overline{\widehat{\sigma}} +\gamma\lambda_2\xi^2 \widehat{\eta}\, \overline{\widehat{v}}\right),
\end{eqnarray}
where 
\begin{equation}\label{46I1I2}
I_1 =\lambda_4\xi^2 -\lambda_2 -\lambda_{9}\quad\hbox{and}\quad I_2 = \lambda_5 \xi^2 -\lambda_2 - \lambda_8.
\end{equation}
To eliminate the terms $Re\left(i\xi\widehat{\theta}\, \overline{\widehat{u}}\right)$, 
$Re\left(i\xi\widehat{y}\, \overline{\widehat{u}}\right)$ and $Re\left(\widehat{\sigma}\, \overline{\widehat{u}}\right)$ from the right hand side of \eqref{g26+}, we put  
\begin{equation*}
I_3 = \gamma +\frac{\vert\gamma\vert}{\gamma}k_4\lambda_0 +\frac{k_4}{\gamma}I_1 ,\quad I_4 = \gamma +\frac{\vert\gamma\vert}{\gamma}k_4\lambda_0 +\frac{k_4}{\gamma}I_2 \quad\hbox{and}\quad I_5 = \gamma +\frac{\vert\gamma\vert}{\gamma}k_4\lambda_0 , 
\end{equation*}
and introduce the functional
\begin{eqnarray}\label{F}
F_1 (\xi,\,t) &=& F_0 (\xi,\,t)+\frac{\vert\gamma\vert}{\gamma}\lambda_0 \,Re (i\xi \widehat{u}\, \overline{\widehat{\eta}})-\frac{1}{\gamma} Re\left(I_1\,\widehat{\eta}\, \overline{\widehat{\theta}} +I_2\,\widehat{\eta}\, \overline{\widehat{y}}\right) \nonumber \\ 
&& +I_5\,Re (i\xi \widehat{v}\, \overline{\widehat{\sigma}})-I_3 Re\left(\widehat{\sigma}\, \overline{\widehat{\phi}}\right)
-I_4 Re\left(\widehat{\sigma}\, \overline{\widehat{z}}\right).
\end{eqnarray}
Multiplying \eqref{equ1}-\eqref{equ3p} by $\lambda_0$, $-\frac{1}{\gamma} I_1$, $-\frac{1}{\gamma} I_2$, 
$I_5$, $I_4$ and $I_3$, respectively, adding the obtained equations, and adding \eqref{g26+}, we arrive at
\begin{eqnarray}\label{Fprime1}
\frac{d}{dt}{F}_1 (\xi,\,t) &=& - \xi^{2}\left(k_2\lambda_1\,\vert\widehat{z}\vert^{2} 
+k_3\lambda_3 \,\vert\widehat{\phi}\vert^{2} +\left(\lambda_5 -\lambda_1 \right)\,\vert\widehat{y}\vert^{2} +\left(\lambda_4 -\lambda_3 \right)\,\vert\widehat{\theta}\vert^{2} +(k_1\lambda_2 -k_1\lambda_4 -k_1 \lambda_5 )\,\vert\widehat{v}\vert^{2}\right)  \nonumber \\ 
&& - \xi^{2}\left((\vert\gamma\vert\lambda_0 -\lambda_2)\vert\widehat{u}\vert^{2}+k_4\,\vert\widehat{\sigma}\vert^{2} \right)+(\vert\gamma\vert\lambda_0 +1 )\xi^{2} \vert\widehat{\eta}\vert^{2} 
+Re\left( i\frac{\vert\gamma\vert}{\gamma}k_5\lambda_0\xi^{2\epsilon_0 +1} \widehat{\eta}\, \overline{\widehat{u}} -ik_5 \xi^{2\epsilon_0 +1}\,\widehat{\eta}\, \overline{\widehat{\sigma}}\right) \nonumber \\
&& +Re\left[\left(\frac{k_1}{\gamma}(I_1 +I_2 )+\left(\gamma\lambda_2 +\gamma +\frac{\vert\gamma\vert}{\gamma}(k_4 -k_1)\lambda_0\right)\xi^2\right)\widehat{\eta}\, \overline{\widehat{v}} +\frac{k_5}{\gamma}I_1 \xi^{2\epsilon_0}\widehat{\eta}\, \overline{\widehat{\theta}}+\frac{k_5}{\gamma}I_2 \xi^{2\epsilon_0}\widehat{\eta}\, \overline{\widehat{y}}\right] \nonumber \\
&& +Re\left[i\left(\gamma\lambda_8 +\frac{k_2}{\gamma}I_2 -I_4\right)\xi\,\widehat{\eta}\, \overline{\widehat{z}} + i\left(\gamma\lambda_{9} +\frac{k_3}{\gamma}I_1 -I_3\right)\xi\widehat{\eta}\, \overline{\widehat{\phi}}\right].
\end{eqnarray}
Let $\lambda$ be a positive constant. We introduce the functionals (${\tilde f}$ is defined in \eqref{22})  
\begin{eqnarray}\label{2FF1}
F (\xi,\,t) =\xi^{2+2\epsilon_0} F_1 (\xi,\,t)\quad\hbox{and}\quad L(\xi,\,t) =\lambda\,\widehat{E}(\xi,\,t) +
\frac{1}{{\tilde f}(\xi)}\,{F}(\xi,\,t).
\end{eqnarray}
Applying Young's inequality for the terms depending on $\widehat{\eta}$ in \eqref{Fprime1}, it follows, for any 
$\varepsilon >0$, that (observe that $I_1 ,\,\cdots, \,I_4$ are of the form $C\xi^2 +{\tilde C}$)
\begin{eqnarray} \label{Fprime2}
\frac{d}{dt}{F} (\xi,\,t) &\leq & - \xi^{4+2\epsilon_0}\left((k_2\lambda_1 -\varepsilon)\,\vert\widehat{z}\vert^{2} 
+(k_3\lambda_3 -\varepsilon)\,\vert\widehat{\phi}\vert^{2} +\left(\lambda_5 -\lambda_1 -\varepsilon\right)\,\vert\widehat{y}\vert^{2} +\left(\lambda_4 -\lambda_3 -\varepsilon\right)\,\vert\widehat{\theta}\vert^{2}\right) \nonumber \\ 
&& - \xi^{4+2\epsilon_0}\left((k_1\lambda_2 -k_1\lambda_4 -k_1 \lambda_5 -\varepsilon) \,\vert\widehat{v}\vert^{2}+ (\vert\gamma\vert\lambda_0 -\lambda_2 -\varepsilon)\vert\widehat{u}\vert^{2}+(k_4-\varepsilon) \,\vert\widehat{\sigma}\vert^{2}
\right)\nonumber \\ 
& & +C_{\varepsilon ,\lambda_0 ,\cdots,\lambda_{9} }{\tilde f}(\xi)\xi^{2\epsilon_0}\vert\widehat{\eta}\vert^{2} .
\end{eqnarray}
We choose $\lambda_{1},\, \lambda_{3} >0$, then we choose $\lambda_0$ such that 
$\lambda_{0} >\frac{1}{\vert\gamma\vert}(\lambda_{1}+\lambda_{3}) $. After, we choose $\lambda_4$ and $\lambda_2$ such that
\begin{equation*}
\lambda_{3} <\lambda_4 <\vert\gamma\vert\lambda_0 -\lambda_1 \quad\hbox{and}\quad \lambda_{1} +\lambda_{4} <\lambda_{2}<\vert\gamma\vert\lambda_0 . 
\end{equation*} 
Finally, we choose $\lambda_5$ and $\varepsilon$ such that
$\lambda_{1} <\lambda_{5}<\lambda_{2}-\lambda_{4}$ and 
\begin{equation*}
0<\varepsilon <\min\left\{\lambda_5 -\lambda_1 ,k_1 (\lambda_2 -\lambda_4 -\lambda_5), \lambda_4 -\lambda_3 , \vert\gamma\vert\lambda_0 -\lambda_2 , k_2\lambda_1 , k_3 \lambda_3 , k_4 \right\}. 
\end{equation*} 
Hence, using the definition \eqref{E21} of $\widehat{E}$, \eqref{Fprime2} leads to, for some positive constant $c_1$,
\begin{eqnarray}\label{g27+} 
\frac{d}{dt}{F}(\xi,\,t) \leq -c_1 \xi^{4+2\epsilon_0}\widehat{E} (\xi,\,t)+ 
C\,{\tilde f}(\xi)\xi^{2\epsilon_0}\vert\widehat{\eta}\vert^{2} .
\end{eqnarray}
Thus, from \eqref{Ep21}, \eqref{2FF1} and \eqref{g27+}, we have 
\begin{eqnarray}\label{g30}
\frac{d}{dt}{L}(\xi,\,t) \leq -c_1 f(\xi)\widehat{E}(\xi,\,t)-\left(k_5\,\lambda - C\right)\xi^{2\epsilon_0} \vert\widehat{\eta}\vert^{2} , 
\end{eqnarray}
where $f$ is defined in \eqref{22}. Moreover, using the definitions of $\widehat{E} ,\,{F} ,\, L$ and ${\tilde f}$, we get, for some $c_{2} >0$ (independent of $\lambda$),
\begin{equation}\label{l1e} 
\vert{L}(\xi,\,t) - \lambda\,\widehat{E}(\xi,\,t)\vert = \frac{1}{{\tilde f}(\xi)}\,\vert{F}(\xi,\,t)\vert \leq C\frac{(1+\vert\xi\vert +\xi^2 )\xi^{2+2\epsilon_0}}{{\tilde f}(\xi)}\leq c_{2}\,\widehat{E}(\xi,\,t).
\end{equation} 
Therefore, for $\lambda$ large enough so that $\lambda >\max\left\{\frac{C}{k_5},\,c_{2}\right\}$, we deduce from \eqref{g30} and \eqref{l1e} that 
\begin{equation}\label{g31}
\frac{d}{dt}{L}(\xi,\,t) + c_1 \,f(\xi)\,\widehat{E}(\xi,\,t) \leq 0
\end{equation}
and
\begin{equation}\label{g32}
c_{3}\,\widehat{E}(\xi,\,t) \leq {L}(\xi,\,t) \leq c_{4}\,\widehat{E}(\xi,\,t),
\end{equation}
where $c_{3} =\lambda - c_{2} >0$ and $c_{4} =\lambda + c_{2} >0$. Consequently, a combination of \eqref{g31} and the second inequality in \eqref{g32} leads to, for $c=\frac{c_1}{c_4}$,  
\begin{equation}\label{g32+}
\frac{d}{dt}{L} (\xi,\,t) + c\,f(\xi)\,{L} (\xi,\,t) \leq 0.
\end{equation}
Finally, by integration \eqref{g32+} with respect to time $t$ and using \eqref{equiv} and \eqref{g32}, \eqref{11} follows with 
${\widetilde c}=\frac{c_4 \alpha_2}{c_3 \alpha_1}$.
\end{proof}
\vskip0,1truecm
\begin{theorem}\label{theorem41}
Assume that $\chi\ne 0$; that is $k_2\ne k_3$. Let $N,\,\ell\in \mathbb{N}^*$ such that $\ell\leq N$, $U_{0}\in H^{N}(\mathbb{R})\cap L^{1}(\mathbb{R})$
and $U$ be a solution of \eqref{e11}. Then, for any $j\in\{0,\,\ldots,\,N-\ell\}$, there exists $c_0 >0$ such that 
\begin{equation}\label{33}
\Vert\partial_{x}^{j}U\Vert_{L^{2}(\mathbb{R})} \leq c_0 \,(1 + t)^{-1/12 - j/6}\,\Vert U_{0} \Vert_{L^{1}(\mathbb{R})} + c_0 \,(1 + t)^{-\ell/2}\,\Vert\partial_{x}^{j+\ell}U_{0} \Vert_{L^{2}(\mathbb{R})},\quad\forall t\in \mathbb{R}_+ 
\end{equation} 
in case \eqref{s21}, and
\begin{equation}\label{33p}
\Vert\partial_{x}^{j}U\Vert_{L^{2}(\mathbb{R})} \leq c_0 \,(1 + t)^{-1/8 - j/4}\,\Vert U_{0} \Vert_{L^{1}(\mathbb{R})} + c_0 \,(1 + t)^{-\ell/2}\,\Vert\partial_{x}^{j+\ell}U_{0} \Vert_{L^{2}(\mathbb{R})},\quad\forall t\in \mathbb{R}_+ 
\end{equation} 
in case \eqref{s22}.  
\end{theorem}
\vskip0,1truecm
\begin{proof} From \eqref{22} we have in case \eqref{s21} (low and high frequences)
\begin{equation}\label{g36}
f(\xi) \geq \left\{
\begin{array}{cc}
\frac{1}{5}\,\xi^{6} & \mbox{if} \quad \vert\xi \vert\leq 1, \\
\\
\frac{1}{5}\,\xi^{-2} & \mbox{if} \quad \vert\xi\vert > 1.
\end{array}
\right.
\end{equation}
Applying Plancherel's theorem and \eqref{11}, we entail
\begin{equation}\label{g37}
\Vert\partial_{x}^{j}U\Vert_{L^{2}(\mathbb{R})}^{2} = \left\Vert\widehat{\partial_{x}^{j} U}(x,\,t)\right \Vert_{L^{2}(\mathbb{R})}^{2}  = \int_{\mathbb{R}}\xi^{2\,j}\,\vert\widehat{U}(\xi,\,t)\vert^{2} d\xi 
\end{equation}
\begin{eqnarray*}
& \leq & \widetilde{c}\int_{\mathbb{R}}\xi^{2\,j}\,e^{-\,c\,f(\xi)\,t}\,\vert\widehat{U}_0 (\xi)\vert^{2} d\xi \\
& \leq & \widetilde{c}\int_{\vert\xi\vert\leq 1}\xi^{2\,j}\,e^{-\,c\,f(\xi)\,t}\,\vert\widehat{U}_0 (\xi)\vert^{2}  d\xi + \widetilde{c}\int_{\vert\xi\vert > 1}\xi^{2\,j}\,e^{-\,c\,f(\xi)\,t}\,\vert\widehat{U}_0 (\xi)\vert^{2}\, d\xi := J_{1} + J_{2}.
\end{eqnarray*}
Using \eqref{g37+} (with $r_1 =2j$, $r_3 =\frac{c}{5}$ and $r_2 =6$) and \eqref{g36}, 
it follows, for the low frequency region,
\begin{equation}\label{g38}
J_{1} \leq C\,\Vert\widehat{U}_{0}\Vert_{L^{\infty}(\mathbb{R})}^{2}\int_{\vert\xi\vert\leq 1}\xi^{2\,j}\,e^{-\,\frac{c}{5}\,t\,\xi^{6}}\ d\xi \leq C\left(1 + t\right)^{-\,\frac{1}{6}(1 + 2\,j)}\Vert {U}_{0}\Vert_{L^{1}(\mathbb{R})}^{2}.
\end{equation}
For the high frequency region, using \eqref{g36}, we observe that
\begin{eqnarray*}
J_{2} & \leq & C\int_{\vert\xi\vert > 1}\vert\xi\vert^{2\,j}\,e^{-\frac{c}{5}\,t\,\xi^{-2}}\,\vert\widehat{U}(\xi,\,0)\vert^{2}\, d\xi \\
& \leq & C\ \sup_{\vert\xi\vert >1}\left\{\vert\xi\vert^{-2\,\ell}\, e^{-\frac{c}{5}\,t\,\vert\xi\vert^{-2}}\right\} \int_{\mathbb{R}}\vert\xi\vert^{2\,(j + \ell)}\,\vert\widehat{U}(\xi,\,0)\vert^{2}\ d\xi,
\end{eqnarray*}
then, using \eqref{ineq} (with $\sigma_1 =2l$, $\sigma_2 =\frac{c}{5}$ and $\sigma_3 =2$),
\begin{eqnarray}\label{g38+}
J_{2} \leq C\left(1 + t\right)^{-\,\ell}\,\Vert\,\partial_{x}^{j + \ell}U_{0}\,\Vert_{L^{2} (\mathbb{R})}^{2} , 
\end{eqnarray}
and so, by combining \eqref{g37}-\eqref{g38+}, we get \eqref{33}.
\vskip0,1truecm
The proof of \eqref{33p} is very similar; we notice only, in case \eqref{s22}, that
\begin{equation*}
f(\xi) \geq \left\{
\begin{array}{cc}
\frac{1}{4}\,\xi^{4} & \mbox{if} \quad \vert\xi \vert\leq 1, \\
\\
\frac{1}{4}\,\xi^{-2} & \mbox{if} \quad \vert\xi\vert > 1.
\end{array}
\right.
\end{equation*}
\end{proof}
\begin{remark}
It is well known, in the literature, that the behavior of the Fourier transform of $U$ in
the low frequency region determines the rate of decay of $U$, while its behavior in the high
frequency rigion imposes a regularity restriction on the initial data known as the regularity-loss property; see 
\cite{4, 8, s4, 6, 10, 7}. The fact that $f$ tends to $0$ when $\xi$ goes to infinity leads to the regularity-loss property in the estimates on $\Vert\partial_{x}^{j}U\Vert_{L^{2}(\mathbb{R})}$ because \eqref{33} and \eqref{33p} with $j=\ell=0$ imply only the boundedness of $\Vert U\Vert_{L^{2}(\mathbb{R})}$. This remark is valid also in case \eqref{tau12} for \eqref{s22}, and in case \eqref{tau12} for \eqref{s21} if \eqref{k123} is not satisfied (see Theorem \ref{theorem441} and Theorem \ref{theorem541} below).
\end{remark} 
\vskip0,1truecm
We finish this subsection by proving that \eqref{g6} is not stable if $(\tau_1 ,\tau_2 ,\tau_3 )=(1,0,0)$ and $\chi=0$. 
\vskip0,1truecm
\begin{theorem}\label{theorem410}
Assume that $\chi=0$; that is $k_2 = k_3$. Then $\vert\widehat{U}(\xi,\,t)\vert$ doesn't converge to zero when time $t$ goes to infinity. 
\end{theorem}
\vskip0,1truecm
\begin{proof} 
We show that, for any $\xi\in \mathbb{R}$, the matrix 
\begin{equation}\label{A012}
A:=-(-\xi^2 A_2 +i\xi A_1 +A_0)
\end{equation}
has at least a pure imaginary eigenvalue; that is   
\begin{equation}\label{AlambdaI}
\forall \xi\in \mathbb{R},\,\,\exists \lambda\in \mathbb{C} :\,\,Re\,(\lambda ) =0,\quad Im\,(\lambda )\ne 0\quad\hbox{and}\quad det(\lambda I -A)=0, 
\end{equation}
where $I$ denotes the identity matrix. From \eqref{e12} with $(\tau_1 ,\tau_2 ,\tau_3 )=(1,0,0)$ and $k_2 =k_3$, we have
\begin{equation*}
\lambda I -A = 
\begin{pmatrix}
\lambda & -i\xi & 0 & -1 & 0 & -1 & 0 & 0 \\
-ik_1 \xi & \lambda & 0 & 0 & 0 & 0 & 0 & i\gamma \xi \\
0 & 0 & \lambda & -i\xi & 0 & 0 & 0 & 0 \\
k_1 & 0 & -ik_2\xi & \lambda & 0 & 0 & 0 & 0 \\
0 & 0 & 0 & 0 & \lambda & -i\xi & 0 & 0 \\
k_1 & 0 & 0 & 0 & -ik_2 \xi & \lambda & 0 & 0 \\
0 & 0 & 0 & 0 & 0 & 0 & \lambda & -i\xi \\
0 & i\gamma\xi & 0 & 0 & 0 & 0 & -ik_4\xi & k_5\xi^{2\epsilon_0} +\lambda
\end{pmatrix}.
\end{equation*}
A direct computaion shows that 
\begin{eqnarray*}
det(\lambda I -A) &=& 2k_1\lambda^2 (\lambda^2 +k_2 \xi^2 )\left[\lambda (\lambda +k_5 \xi^{2\epsilon_0})+(k_4 +\gamma^2 )\xi^2 \right] +k_4 \xi^2 (\lambda^2 +k_1 \xi^2 )\left(\lambda^2 + k_2\xi^2 \right)^2\\
&& +\lambda (\lambda^2 +k_2 \xi^2 )^2 \left[\lambda^2 (\lambda +k_5 \xi^{2\epsilon_0})+\gamma^2 \lambda \xi^2 +k_1 \xi^2 (\lambda +k_5 \xi^{2\epsilon_0})\right] . 
\end{eqnarray*}
It is clear that, if $\xi \ne 0$, then $\lambda =i {\sqrt{k_2}}\xi$ is a pure imaginary eigenvalue of $A$.
If $\xi = 0$, then $\lambda =i{\sqrt{2k_1}}$ is a pure imaginary eigenvalue of $A$. Consequently, according to \eqref{e14} and \eqref{A012} (see  \cite{tesc}), the solution of \eqref{g6} doesn't converge to zero when time $t$ goes to infinity.
\end{proof} 

\subsection{Case 2: $(\tau_1 ,\tau_2 ,\tau_3)=(0,1,0)$}

We present, first, our exponential stability result for \eqref{g6}, where the proof is similar to the one of Lemma \ref{lemma32}. 
\vskip0,1truecm
\begin{lemma}\label{lemma432}
Let $\widehat{U}$ be a solution of \eqref{g6}. Then, there exist $c,\,\widetilde{c}>0$ such that \eqref{11} is satisfied with  
\begin{align}\label{422}
f(\xi) =\frac{\xi^{4+2\epsilon_0 }}{{\tilde f}(\xi)}\quad\hbox{and}\quad {\tilde f}(\xi) = \begin{cases}
1 +\xi^{2} +\xi^{4} +\xi^{6}\quad &\hbox{for \eqref{s21} and \eqref{s22} under \eqref{k123}},\\
1 +\xi^{2} +\xi^4 +\xi^6 +\xi^{8} +\xi^{10} \quad &\hbox{for \eqref{s21} without \eqref{k123}},\\
1 +\xi^{2} +\xi^4 +\xi^6 +\xi^{8} \quad &\hbox{for \eqref{s22} without \eqref{k123}.}
\end{cases}
\end{align}
\end{lemma}
\vskip0,1truecm
\begin{proof}
Multiplying \eqref{e101}$_{4}$ and \eqref{e101}$_{8}$ by $i\frac{\vert\gamma\vert}{\gamma}\,\xi\,\overline{\widehat{\eta}}$ and $-i\frac{\vert\gamma\vert}{\gamma}\,\xi\,\overline{\widehat{y}}$, respectively, adding the resulting equations, taking the real part and using \eqref{hd1}, we get
\begin{eqnarray}\label{equp12}
\frac{d}{dt}Re\left(i\frac{\vert\gamma\vert}{\gamma}\,\xi\,\widehat{y}\, \overline{\widehat{\eta}}\right) &=& \vert\gamma\vert\xi^{2}\left(\vert\widehat{\eta}\vert^{2} - \vert\widehat{y}\vert^{2}\right) + \frac{\vert\gamma\vert}{\gamma}k_4\xi^2\,Re\left(\widehat{\sigma}\,\overline{\widehat{y}}\right) 
 -\frac{\vert\gamma\vert}{\gamma}k_1 \,Re\left(i\xi\widehat{v}\, \overline{\widehat{\eta}}\right) \nonumber \\ 
&& -\frac{\vert\gamma\vert}{\gamma}k_2 \xi^2\,Re\left(\widehat{\eta}\, \overline{\widehat{z}}\right) +\frac{\vert\gamma\vert}{\gamma}k_5 \xi^{2\epsilon_0}\,Re\left(i\xi\widehat{\eta}\, \overline{\widehat{y}}\right).
\end{eqnarray}
Multiplying \eqref{e101}$_{1}$ and \eqref{e101}$_{7}$ by $\xi^2\overline{\widehat{\sigma}}$ and $\xi^2\overline{\widehat{v}}$, respectively, adding the resulting equations, taking the real part and using \eqref{hd1}, we find
\begin{equation}\label{equp22}
\frac{d}{dt}Re\left(\xi^2\widehat{v}\, \overline{\widehat{\sigma}}\right) = \xi^2 Re\left(\widehat{\sigma}\,\overline{\widehat{y}}\right)+\xi^2\,Re\left(\widehat{\sigma}\, \overline{\widehat{\theta}}\right) +\xi^2 \,Re\left(i\xi\widehat{u}\, \overline{\widehat{\sigma}}\right) +\xi^2\,Re\left(i\xi\widehat{\eta}\, \overline{\widehat{v}}\right).
\end{equation}
Also, multiplying \eqref{e101}$_{3}$ and \eqref{e101}$_{7}$ by $i\xi\overline{\widehat{\sigma}}$ and $-i\xi\overline{\widehat{z}}$, respectively, adding the resulting equations, taking the real part and using \eqref{hd1}, we obtain
\begin{equation}\label{equp32}
\frac{d}{dt}Re\left(i\xi\widehat{z}\, \overline{\widehat{\sigma}}\right) = -\xi^2 Re\left(\widehat{\sigma}\,\overline{\widehat{y}}\right)+\xi^2\,Re\left(\widehat{\eta}\, \overline{\widehat{z}}\right).
\end{equation}
Multiplying \eqref{e101}$_{2}$ and \eqref{e101}$_{8}$ by $\overline{\widehat{\eta}}$ and $\overline{\widehat{u}}$, respectively, adding the resulting equations, taking the real part and using \eqref{hd1}, we infer that
\begin{equation}\label{equp42}
\frac{d}{dt}Re\left(\widehat{u}\, \overline{\widehat{\eta}}\right) = -\gamma Re\left(i\xi\widehat{y}\,\overline{\widehat{u}}\right)+k_1 Re\left(i\xi\widehat{v}\, \overline{\widehat{\eta}}\right) +k_4 Re\left(i\xi\widehat{\sigma}\, \overline{\widehat{u}}\right) -k_5 \xi^{2\epsilon_0}Re\left(\widehat{\eta}\, \overline{\widehat{u}}\right).
\end{equation}
Multiplying \eqref{e101}$_{5}$ and \eqref{e101}$_{7}$ by $i\xi\overline{\widehat{\sigma}}$ and $-i\xi\overline{\widehat{\phi}}$, respectively, adding the resulting equations, taking the real part and using \eqref{hd1}, we infer that
\begin{equation}\label{equp52}
\frac{d}{dt}Re\left(i\xi\widehat{\phi}\, \overline{\widehat{\sigma}}\right) = -\xi^2 Re\left(\widehat{\sigma}\,\overline{\widehat{\theta}}\right)+\xi^2 Re\left(\widehat{\eta}\, \overline{\widehat{\phi}}\right) .
\end{equation}
Finally, multiplying \eqref{e101}$_{6}$ and \eqref{e101}$_{8}$ by $-i\xi\overline{\widehat{\eta}}$ and $i\xi\overline{\widehat{\theta}}$, respectively, adding the resulting equations, taking the real part and using \eqref{hd1}, it follows that
\begin{equation}\label{equp62}
\frac{d}{dt}Re\left(i\xi\widehat{\eta}\, \overline{\widehat{\theta}}\right) = \gamma\xi^2 Re\left(\widehat{y}\,\overline{\widehat{\theta}}\right)-k_4 \xi^2 Re\left(\widehat{\sigma}\, \overline{\widehat{\theta}}\right) +k_3\xi^2 Re\left(\widehat{\eta}\, \overline{\widehat{\phi}}\right) -k_5 \xi^{2\epsilon_0}Re\left(i\xi\widehat{\eta}\, \overline{\widehat{\theta}}\right)+k_1 Re\left(i\xi\widehat{v}\, \overline{\widehat{\eta}}\right) .
\end{equation}
Let $\lambda_{0} ,\,\cdots ,\,\lambda_{5}$ be positive constants, and let  
\begin{equation*}
\lambda_{6} =\frac{k_2}{k_3}\left[\left(\frac{k_3}{k_1} -1\right)\lambda_4 \xi^2 -\lambda_{2}-\lambda_{3}\right],\quad \lambda_{7} =-\frac{k_3}{k_2}\lambda_{6}, \quad \lambda_{8} =-\frac{k_2}{k_1}\lambda_{5} \xi^2 +\lambda_{6} -\lambda_{1} \quad\hbox{and}\quad
\lambda_{9} = \lambda_{4} \xi^2 +\lambda_2 .
\end{equation*}
We define the functional ${F}_0$ as follows:
\begin{eqnarray}\label{F02}
{F}_0 (\xi,\,t) & = & Re\left[i\,\xi\,\left(\lambda_1\,\widehat{y}\, \overline{\widehat{z}}- \lambda_2\,\widehat{u}\, \overline{\widehat{v}} +\lambda_3\widehat{\theta}\, \overline{\widehat{\phi}} +\widehat{\eta}\, \overline{\widehat{\sigma}} +\lambda_6\widehat{z}\, \overline{\widehat{\theta}}+\lambda_7\widehat{\phi}\,\overline{\widehat{y}}\right)\right] \nonumber \\ 
& & +Re\left(-\lambda_4\xi^2 \,\widehat{\theta}\, \overline{\widehat{v}} +\lambda_5\xi^2 \,\widehat{y}\, \overline{\widehat{v}}-\lambda_8 \widehat{u}\, \overline{\widehat{z}}
-\lambda_{9} \widehat{u}\, \overline{\widehat{\phi}}\right).
\end{eqnarray}
By multiplying \eqref{eq31}-\eqref{eq312} by $\lambda_1 ,\,-\lambda_2 ,\,\lambda_3 ,\,\lambda_4$, $-\lambda_5 ,\,1$ and 
$\lambda_6 ,\,\cdots,\,\lambda_{9}$, respectively, and adding the resulting equations, we find
\begin{eqnarray}\label{g26+2}
\frac{d}{dt}{F}_0 (\xi,\,t) &=& - \xi^{2}\left(k_3\lambda_3 \,\vert\widehat{\phi}\vert^{2} +\lambda_2\,\vert\widehat{u}\vert^{2} +\left(\lambda_4 -\lambda_3 \right)\,\vert\widehat{\theta}\vert^{2} +(k_1\lambda_5 -k_1\lambda_4 -k_1 \lambda_2 )\,\vert\widehat{v}\vert^{2}\right)  \nonumber \\ 
&& - \xi^{2}\left(k_2\lambda_1\,\vert\widehat{z}\vert^{2} 
+k_4\,\vert\widehat{\sigma}\vert^{2} \right)
+I_1 Re (i\xi\widehat{y}\, \overline{\widehat{u}})+I_2 \xi^2 Re (\widehat{y}\, \overline{\widehat{\theta}}) +
\gamma\xi^{2} Re (\widehat{\sigma}\, \overline{\widehat{y}}) \nonumber \\ 
&& +\xi^{2} \left((\lambda_1 +\lambda_5)\vert\widehat{y}\vert^{2}
+\vert\widehat{\eta}\vert^{2}\right) +Re\left( \gamma\lambda_1\xi^2\,\widehat{\eta}\, \overline{\widehat{z}} - \gamma\lambda_{7}\xi^2\widehat{\eta}\, \overline{\widehat{\phi}} -ik_5 \xi^{2\epsilon_0 +1}\,\widehat{\eta}\, \overline{\widehat{\sigma}} -i\gamma\lambda_5\xi^3 \widehat{\eta}\, \overline{\widehat{v}}\right)
\end{eqnarray}
(thanks to the choices of $\lambda_6 ,\,\cdots ,\,\lambda_9$, $Re \left(i\xi\left(\widehat{v}\, \overline{\widehat{z}} +\widehat{\theta}\, \overline{\widehat{u}}+\widehat{v}\, \overline{\widehat{\phi}}\right)\right)$ and $Re\left(\widehat{z}\, \overline{\widehat{\phi}}\right)$ disappear), where 
\begin{equation*}
I_1 =-\lambda_5\xi^2 +\lambda_2 -\lambda_{8}\quad\hbox{and}\quad I_2 = \lambda_5 -\lambda_4 - \lambda_6 - \lambda_7 .
\end{equation*}
We put  
\begin{equation*}
I_3 =\left(\frac{\vert\gamma\vert}{\gamma}k_4\lambda_0 +\gamma\right)\xi^2 +\frac{k_4}{\gamma}I_1\quad\hbox{and}\quad I_4 = \frac{k_4}{\gamma}(I_2\xi^2 +I_1 ), 
\end{equation*}
and introduce the functional
\begin{eqnarray}\label{F}
F_1 (\xi,\,t) &=& \xi^2 F_0 (\xi,\,t)+\frac{\vert\gamma\vert}{\gamma}\lambda_0 \xi^2\,Re (i\xi \widehat{y}\, \overline{\widehat{\eta}})+\frac{k_4}{\gamma} I_1 \xi^2 Re\left(\widehat{v}\, \overline{\widehat{\sigma}}\right) 
+I_3 Re\left(i\xi\widehat{z}\, \overline{\widehat{\sigma}}\right)\nonumber \\ 
&& +\frac{1}{\gamma} I_1 \xi^2\,Re (\widehat{u}\, \overline{\widehat{\eta}})+I_4 Re\left(i\xi\widehat{\phi}\, \overline{\widehat{\sigma}}\right)-\frac{1}{\gamma} I_2 \xi^2 Re\left(i\xi\widehat{\eta}\, \overline{\widehat{\theta}}\right).
\end{eqnarray}
Multiplying \eqref{equp12}-\eqref{equp62} and \eqref{g26+2} by $\lambda_0 \xi^2$, $\frac{k_4}{\gamma} I_1$, $I_3$, 
$\frac{1}{\gamma} I_1 \xi^2$, $I_4$, $-\frac{1}{\gamma} I_2 \xi^2$ and $\xi^2$, respectively, and adding the obtained expressions, we arrive at ($Re \left(i\xi\widehat{u}\, \overline{\widehat{\sigma}}\right)$ and $Re \left(\widehat{\sigma}\, \overline{\widehat{y}} +\widehat{\sigma}\, \overline{\widehat{\theta}}\right)$ disappear according to the definition of $I_3$ and $I_4$)
\begin{eqnarray}\label{FFprime1}
\frac{d}{dt}{F}_1 (\xi,\,t) &=& - \xi^{4}\left(k_2\lambda_1\,\vert\widehat{z}\vert^{2} 
+k_3\lambda_3 \,\vert\widehat{\phi}\vert^{2} +\lambda_2\,\vert\widehat{u}\vert^{2} +\left(\lambda_4 -\lambda_3 \right)\,\vert\widehat{\theta}\vert^{2} +(k_1\lambda_5 -k_1\lambda_4 -k_1 \lambda_2 )\,\vert\widehat{v}\vert^{2}\right)  \nonumber \\ 
&& - \xi^{4}\left((\vert\gamma\vert\lambda_0 -\lambda_1 -\lambda_5)\vert\widehat{y}\vert^{2}+k_4\,\vert\widehat{\sigma}\vert^{2} \right)+(\vert\gamma\vert\lambda_0 +1 )\xi^{4} \vert\widehat{\eta}\vert^{2} \nonumber \\
&& +\xi^2 Re\left[\left(iI_5 \overline{\widehat{v}}+I_6 \overline{\widehat{z}}+I_7 \overline{\widehat{\phi}}-ik_5 \xi^{2\epsilon_0 +1}\overline{\widehat{\sigma}}+i\frac{\vert\gamma\vert}{\gamma}k_5 \lambda_0 \xi^{2\epsilon_0 +1}\overline{\widehat{y}}
+i\frac{k_5}{\gamma} \xi^{2\epsilon_0 +1}I_2 \overline{\widehat{\theta}} -\frac{k_5}{\gamma} \xi^{2\epsilon_0 }I_1 \overline{\widehat{u}}\right) \widehat{\eta}\right],
\end{eqnarray}
where
\begin{equation*}
I_5 =-\gamma\lambda_5\xi^3 +\left(\frac{\vert\gamma\vert}{\gamma}k_1\lambda_0 +\frac{k_4 -k_1}{\gamma} I_1 +\frac{k_1}{\gamma} I_2\right)\xi, 
\end{equation*}
\begin{equation*}
I_6 =\left(-\frac{\vert\gamma\vert}{\gamma}k_2\lambda_0 +\gamma\lambda_1 \right)\xi^2 +I_3\quad\hbox{and}\quad I_7 = -\left(\frac{k_3}{\gamma} I_2 +\gamma \lambda_7\right)\xi^2 + I_4 . 
\end{equation*}
Observe that, by definition,
\begin{equation}\label{I12}
\vert I_1\vert \leq \begin{cases}
C\quad &\hbox{if \eqref{k123} holds},\\
C(1+\xi^{2} )\quad &\hbox{if not},
\end{cases} \quad \vert I_2\vert \leq \begin{cases}
C\quad &\hbox{if \eqref{k123} holds},\\
C(1+\xi^{2} )\quad &\hbox{if not},
\end{cases} 
\end{equation}
\begin{equation}\label{I57}
\vert I_5\vert \leq C(\vert\xi\vert +\vert\xi\vert^3),\quad \vert I_6\vert \leq C(1+\xi^2) \quad\hbox{and}\quad \vert I_7\vert \leq \begin{cases}
C(1+\xi^{2})\quad &\hbox{if \eqref{k123} holds},\\
C(1+\xi^{2} +\xi^{4})\quad &\hbox{if not}.
\end{cases} 
\end{equation}
Then, applying Young's inequality, it follows, for any $\varepsilon >0$, that
\begin{equation*}
\xi^2 Re\left[\left(iI_5 \overline{\widehat{v}}+I_6 \overline{\widehat{z}}+I_7 \overline{\widehat{\phi}}-ik_5 \xi^{2\epsilon_0 +1}\overline{\widehat{\sigma}}+i\frac{\vert\gamma\vert}{\gamma}k_5 \lambda_0 \xi^{2\epsilon_0 +1}\overline{\widehat{y}}
+i\frac{k_5}{\gamma} \xi^{2\epsilon_0 +1}I_2 \overline{\widehat{\theta}} -\frac{k_5}{\gamma} \xi^{2\epsilon_0 }I_1 \overline{\widehat{u}}\right) \widehat{\eta}\right]
\end{equation*}
\begin{eqnarray}\label{I17}
&\leq& \varepsilon\xi^4 \left(\vert\widehat{z}\vert^{2} +\vert\widehat{\phi}\vert^{2} +\vert\widehat{u}\vert^{2} +\vert\widehat{\theta}\vert^{2} +\vert\widehat{v}\vert^{2}+\vert\widehat{\sigma}\vert^{2} +\vert\widehat{y}\vert^{2}\right) \nonumber \\
&& + C_{\varepsilon} \left(\xi^{4\epsilon_0}\vert I_1\vert^2 +\xi^{4\epsilon_0 +2}\vert I_2\vert^2 +\vert I_5\vert^2 +\vert I_6\vert^2 +\vert I_7\vert^2 +\xi^{4\epsilon_0 +2}\right)\vert\widehat{\eta}\vert^{2} \nonumber \\
&\leq& \varepsilon\xi^4 \left(\vert\widehat{z}\vert^{2} +\vert\widehat{\phi}\vert^{2} +\vert\widehat{u}\vert^{2} +\vert\widehat{\theta}\vert^{2} +\vert\widehat{v}\vert^{2}+\vert\widehat{\sigma}\vert^{2} +\vert\widehat{y}\vert^{2}\right)+
C_{\varepsilon,\lambda_0 ,\cdots ,\lambda_9 }{\tilde f}(\xi)\vert\widehat{\eta}\vert^{2}.
\end{eqnarray}
By combining \eqref{FFprime1} and \eqref{I17}, we find 
\begin{eqnarray} \label{FFprime2}
\frac{d}{dt}{F}_1 (\xi,\,t) &\leq& - \xi^{4}\left((k_2\lambda_1 -\varepsilon)\,\vert\widehat{z}\vert^{2} 
+(k_3\lambda_3 -\varepsilon)\,\vert\widehat{\phi}\vert^{2} +(\lambda_2 -\varepsilon)\,\vert\widehat{u}\vert^{2} +\left(\lambda_4 -\lambda_3 -\varepsilon\right)\,\vert\widehat{\theta}\vert^{2} \right)  \nonumber \\ 
&& - \xi^{4}\left((k_1\lambda_5 -k_1\lambda_4 -k_1 \lambda_2 -\varepsilon )\,\vert\widehat{v}\vert^{2} +(\vert\gamma\vert\lambda_0 -\lambda_1 -\lambda_5 -\varepsilon ) \vert\widehat{y}\vert^{2} +(k_4 -\varepsilon)\vert\widehat{\sigma} \vert^{2}\right) \nonumber \\
&& +C_{\varepsilon,\lambda_0 ,\cdots ,\lambda_9 } {\tilde f}(\xi)\vert\widehat{\eta}\vert^{2}.
\end{eqnarray}
Let $\lambda$ be a positive constant. We introduce the functionals  
\begin{eqnarray}\label{FF1}
F (\xi,\,t) =\xi^{2\epsilon_0} F_1 (\xi,\,t)\quad\hbox{and}\quad L(\xi,\,t) =\lambda\,\widehat{E}(\xi,\,t) +
\frac{1}{{\tilde f}(\xi)}\,{F}(\xi,\,t).
\end{eqnarray}
We choose $0<\lambda_{1}$, $0<\lambda_{3}<\lambda_{4} <\lambda_5$,  $0<\lambda_{2} <\lambda_{5} -\lambda_{4}$, 
$\lambda_{0} >\frac{1}{\vert\gamma\vert}(\lambda_{1} +\lambda_{5})$ and 
\begin{equation*}
0<\varepsilon <\min\left\{k_2 \lambda_1 , k_3 \lambda_3 ,\lambda_2 ,\lambda_4 -\lambda_3 ,
k_1\lambda_5 -k_1\lambda_4 -k_1 \lambda_2 ,\vert\gamma\vert\lambda_0 -\lambda_1 -\lambda_5 , k_4 \right\}, 
\end{equation*}
and use the definition of $\widehat{E}$, we deduce from \eqref{FFprime2} and \eqref{FF1}, for some positive constant $c_1$, that
\begin{equation}\label{4g27+} 
\frac{d}{dt}{F}(\xi,\,t) \leq -c_1 \xi^{4+2\epsilon_0}\widehat{E} (\xi,\,t)+ C{\tilde f}(\xi)\xi^{2\epsilon_0} \vert\widehat{\eta}\vert^{2} .
\end{equation}
Then, from \eqref{Ep21}, \eqref{FF1} and \eqref{4g27+}, we infer that
\begin{eqnarray}\label{4g30}
\frac{d}{dt}{L}(\xi,\,t) \leq -c_1 f(\xi)\widehat{E}(\xi,\,t)-\left(k_5\,\lambda - C\right)\xi^{2\epsilon_0} \vert\widehat{\eta}\vert^{2} . 
\end{eqnarray}
On the other hand, the definitions of $\widehat{E}$, $F$ and $L$ imply that there exists $c_2>0$ (not depending on $\lambda$) such that, for $d_0 =0$ if \eqref{k123} holds, and $d_0 =5$ if not,
\begin{equation*}
\left\vert {L}(\xi,\,t)-\lambda\widehat{E}(\xi,\,t)\right\vert \leq 
c_2 \frac{\xi^{2\epsilon_0}(1+\vert\xi\vert +\xi^{2} +\vert\xi\vert^3 +\xi^{4} +\vert\xi\vert^{d_0})}{{\tilde f}(\xi)}\widehat{E}(\xi,\,t) \leq 6c_2 \widehat{E}(\xi,\,t). 
\end{equation*}
So, we choose $\lambda >\max\left\{\frac{C}{k_5},6c_2\right\}$, we get \eqref{g31} and \eqref{g32} with $c_3 =\lambda -6c_2 >0$ and
$c_4 =\lambda +6c_2 >0$. The proof can be ended as for Lemma \ref{lemma32}.
\end{proof}
\vskip0,1truecm
\begin{theorem}\label{theorem441}
Let $N,\,\ell\in \mathbb{N}^*$ such that $\ell\leq N$, $U_{0}\in H^{N}(\mathbb{R})\cap L^{1}(\mathbb{R})$
and $U$ be the solution of \eqref{e11}. Then, for any $j\in\{0,\,\ldots,\,N-\ell\}$, there exist $c_0 ,\,{\tilde c}_0 >0$ such that, for any 
$t\in \mathbb{R}_+$, 
\vskip0,1truecm
{\bf Case \eqref{s21}}:
\begin{equation}\label{33case20}
\Vert\partial_{x}^{j}U\Vert_{L^{2}(\mathbb{R})} \leq c_0 \,(1 + t)^{-1/12 - j/6}\,\Vert U_{0} \Vert_{L^{1}(\mathbb{R})} + c_0 e^{-{\tilde c}_0 t}\,\Vert\partial_{x}^{j+\ell}U_{0} \Vert_{L^{2}(\mathbb{R})} \quad\hbox{if}\,\,k_1 =k_2 =k_3 , 
\end{equation} 
and 
\begin{equation}\label{33case2}
\Vert\partial_{x}^{j}U\Vert_{L^{2}(\mathbb{R})} \leq c_0 \,(1 + t)^{-1/12 - j/6}\,\Vert U_{0} \Vert_{L^{1}(\mathbb{R})} + c_0 \,(1 + t)^{-\ell/4}\,\Vert\partial_{x}^{j+\ell}U_{0} \Vert_{L^{2}(\mathbb{R})} \quad\hbox{if not}. 
\end{equation} 
\vskip0,1truecm
{\bf Case \eqref{s22}}:
\begin{equation}\label{33case20+}
\Vert\partial_{x}^{j}U\Vert_{L^{2}(\mathbb{R})} \leq c_0 \,(1 + t)^{-1/8 - j/4}\,\Vert U_{0} \Vert_{L^{1}(\mathbb{R})} + c_0 \,(1 + t)^{-\ell/2}\,\Vert\partial_{x}^{j+\ell}U_{0} \Vert_{L^{2}(\mathbb{R})} \quad\hbox{if}\,\,k_1 =k_2 =k_3 , 
\end{equation} 
and 
\begin{equation}\label{33case2+}
\Vert\partial_{x}^{j}U\Vert_{L^{2}(\mathbb{R})} \leq c_0 \,(1 + t)^{-1/8 - j/4}\,\Vert U_{0} \Vert_{L^{1}(\mathbb{R})} + c_0 \,(1 + t)^{-\ell/4}\,\Vert\partial_{x}^{j+\ell}U_{0} \Vert_{L^{2}(\mathbb{R})} \quad\hbox{if not}. 
\end{equation}
\end{theorem}
\vskip0,1truecm
\begin{proof} For \eqref{s21}, we have, from \eqref{422} (low and high frequences),
\begin{equation}
\label{4g3600}f(\xi) \geq \left\{
\begin{array}{cc}
\frac{1}{4}\,\xi^{6} & \mbox{if} \quad \vert\xi \vert\leq 1, \\
\\
\frac{1}{4} & \mbox{if} \quad \vert\xi\vert > 1
\end{array}
\right. \quad\hbox{if}\,\,k_1 =k_2 =k_3 , 
\end{equation}
and
\begin{equation}
\label{4g36}f(\xi) \geq \left\{
\begin{array}{cc}
\frac{1}{6}\,\xi^{6} & \mbox{if} \quad \vert\xi \vert\leq 1, \\
\\
\frac{1}{6}\,\xi^{-4} & \mbox{if} \quad \vert\xi\vert > 1
\end{array}
\right. \quad\hbox{if not}.
\end{equation}
The proof of \eqref{33case2} is identical to the one of Theorem \ref{theorem41} by using \eqref{4g36} and applying \eqref{g37+} (with $r_1 =2j$, 
$r_2 =\frac{c}{6}$ and $r_3 =6$) and \eqref{ineq} (with $\sigma_1 =2l$, $\sigma_2 =\frac{c}{6}$ and $\sigma_3 =4$). To get \eqref{33case20}, noticing that the low frequencies can be treated as for \eqref{33case2}. For the high frequencies, we have just to remark that \eqref{4g3600} implies that
\begin{eqnarray*}
\int_{\vert\xi\vert > 1}\vert\xi\vert^{2\,j}\,e^{-cf(\xi)t}\,\vert\widehat{U}(\xi,\,0)\vert^{2}\ d\xi& \leq & \int_{\vert\xi\vert > 1}\vert\xi\vert^{2\,j}\,e^{-\frac{c}{4}\,t}\,\vert\widehat{U}(\xi,\,0)\vert^{2}\ d\xi  \nonumber \\
& \leq & \sup_{\vert\xi\vert >1}\left\{\vert\xi\vert^{-2\,\ell}\, e^{-\frac{c}{4}\,t}\right\} \int_{\mathbb{R}}\vert\xi\vert^{2\,(j + \ell)}\,\vert\widehat{U}(\xi,\,0)\vert^{2}\ d\xi  \nonumber \\
& \leq & e^{-\frac{c}{4}\,t}\,\Vert\,\partial_{x}^{j + \ell}U_{0}\,\Vert_{L^{2} (\mathbb{R})}^{2} ,
\end{eqnarray*}
so \eqref{33case20} holds true with ${\tilde c}_0 =\frac{c}{8}$. The proof of \eqref{33case20+} and \eqref{33case2+} is identical to the one of
\eqref{33case2} by remarking, for \eqref{s22}, that
\begin{equation*}
f(\xi) \geq \left\{
\begin{array}{cc}
\frac{1}{4}\,\xi^{4} & \mbox{if} \quad \vert\xi \vert\leq 1, \\
\\
\frac{1}{4}\,\xi^{-2}& \mbox{if} \quad \vert\xi\vert > 1
\end{array}
\right. \quad\hbox{if}\,\,k_1 =k_2 =k_3 , 
\end{equation*}
and
\begin{equation*}
f(\xi) \geq \left\{
\begin{array}{cc}
\frac{1}{5}\,\xi^{4} & \mbox{if} \quad \vert\xi \vert\leq 1, \\
\\
\frac{1}{5}\,\xi^{-4} & \mbox{if} \quad \vert\xi\vert > 1
\end{array}
\right.\quad\hbox{if not}.
\end{equation*} 
\end{proof}
\vskip0,1truecm
\begin{remark}
In case \eqref{s21} under \eqref{k123}, the fact that $f$ tends to $1$ when $\xi$ goes to infinity allows to avoid the regularity-loss property in the estimate \eqref{33case20} on $\Vert\partial_{x}^{j}U\Vert_{L^{2}(\mathbb{R})}$ because one can take $j=\ell=0$, and the stability of \eqref{e11} is sill satisfied with a decay estimate depending only on $\Vert U_0\Vert_{L^{1}(\mathbb{R})}$ and $\Vert U_0\Vert_{L^{2}(\mathbb{R})}$. This remark is valid also for \eqref{s21} in case $(\tau_1 ,\tau_2 ,\tau_3)=(0,0,1)$ under \eqref{k123} (see Theorem \ref{theorem541} below).
\end{remark} 

\subsection{Case 3: $(\tau_1 ,\tau_2 ,\tau_3)=(0,0,1)$}

In this case, we prove the same stability results for \eqref{g6} and \eqref{e11} that given in the previous subsection, and moreover, the proofs are very similar. 
\vskip0,1truecm
\begin{lemma}\label{lemma532}
The result of Lemma \ref{lemma432} holds true also when $(\tau_1 ,\tau_2 ,\tau_3)=(0,0,1)$. 
\end{lemma}
\vskip0,1truecm
\begin{proof}
Multiplying \eqref{e101}$_{6}$ and \eqref{e101}$_{8}$ by $i\frac{\vert\gamma\vert}{\gamma}\,\xi\,\overline{\widehat{\eta}}$ and $-i\frac{\vert\gamma\vert}{\gamma}\,\xi\,\overline{\widehat{\theta}}$, respectively, adding the resulting equations, taking the real part and using \eqref{hd1}, we get
\begin{eqnarray}\label{equp123}
\frac{d}{dt}Re\left(i\frac{\vert\gamma\vert}{\gamma}\,\xi\,\widehat{\theta}\, \overline{\widehat{\eta}}\right) &=& \vert\gamma\vert\xi^{2}\left(\vert\widehat{\eta}\vert^{2} - \vert\widehat{\theta}\vert^{2}\right) + \frac{\vert\gamma\vert}{\gamma}k_4\xi^2\,Re\left(\widehat{\sigma}\,\overline{\widehat{\theta}}\right) 
 -\frac{\vert\gamma\vert}{\gamma}k_1 \,Re\left(i\xi\widehat{v}\, \overline{\widehat{\eta}}\right) \nonumber \\ 
&& -\frac{\vert\gamma\vert}{\gamma}k_3 \xi^2\,Re\left(\widehat{\eta}\, \overline{\widehat{\phi}}\right) +\frac{\vert\gamma\vert}{\gamma}k_5 \xi^{2\epsilon_0}\,Re\left(i\xi\widehat{\eta}\, \overline{\widehat{\theta}}\right).
\end{eqnarray}
Also, multiplying \eqref{e101}$_{2}$ and \eqref{e101}$_{8}$ by $\overline{\widehat{\eta}}$ and $\overline{\widehat{u}}$, respectively, adding the resulting equations, taking the real part and using \eqref{hd1}, we infer that
\begin{equation}\label{equp423}
\frac{d}{dt}Re\left(\widehat{u}\, \overline{\widehat{\eta}}\right) = -\gamma Re\left(i\xi\widehat{\theta}\,\overline{\widehat{u}}\right)+k_1 Re\left(i\xi\widehat{v}\, \overline{\widehat{\eta}}\right) +k_4 Re\left(i\xi\widehat{\sigma}\, \overline{\widehat{u}}\right) -k_5 \xi^{2\epsilon_0}Re\left(\widehat{\eta}\, \overline{\widehat{u}}\right).
\end{equation}
Finally, multiplying \eqref{e101}$_{4}$ and \eqref{e101}$_{8}$ by $-i\xi\overline{\widehat{\eta}}$ and $i\xi\overline{\widehat{y}}$, respectively, adding the resulting equations, taking the real part and using \eqref{hd1}, it follows that
\begin{equation}\label{equp623}
\frac{d}{dt}Re\left(i\xi\widehat{\eta}\, \overline{\widehat{y}}\right) = \gamma\xi^2 Re\left(\widehat{y}\,\overline{\widehat{\theta}}\right)-k_4 \xi^2 Re\left(\widehat{\sigma}\, \overline{\widehat{y}}\right) +k_2\xi^2 Re\left(\widehat{\eta}\, \overline{\widehat{z}}\right) -k_5 \xi^{2\epsilon_0}Re\left(i\xi\widehat{\eta}\, \overline{\widehat{y}}\right)+k_1 Re\left(i\xi\widehat{v}\, \overline{\widehat{\eta}}\right) .
\end{equation}
Let $\lambda_{0} ,\,\cdots ,\,\lambda_{5}$ be positive constants, and let  
\begin{equation*}
\lambda_{6} =\frac{k_2}{k_3}\left[\left(1-\frac{k_3}{k_1}\right)\lambda_4 \xi^2 -\lambda_{2}-\lambda_{3}\right],\quad \lambda_{7} =-\frac{k_3}{k_2}\lambda_{6}, \quad \lambda_{8} =\frac{k_2}{k_1}\lambda_{5} \xi^2 +\lambda_{6} -\lambda_{1} \quad\hbox{and}\quad
\lambda_{9} = -\lambda_{4} \xi^2 +\lambda_2 .
\end{equation*}
We define the functional ${F}_0$ as follows:
\begin{eqnarray}\label{F023}
{F}_0 (\xi,\,t) & = & Re\left[i\,\xi\,\left(\lambda_1\,\widehat{y}\, \overline{\widehat{z}}- \lambda_2\,\widehat{u}\, \overline{\widehat{v}} +\lambda_3\widehat{\theta}\, \overline{\widehat{\phi}} +\widehat{\eta}\, \overline{\widehat{\sigma}} +\lambda_6\widehat{z}\, \overline{\widehat{\theta}}+\lambda_7\widehat{\phi}\,\overline{\widehat{y}}\right)\right] \nonumber \\ 
& & +Re\left(\lambda_4\xi^2 \,\widehat{\theta}\, \overline{\widehat{v}} -\lambda_5\xi^2 \,\widehat{y}\, \overline{\widehat{v}}-\lambda_8 \widehat{u}\, \overline{\widehat{z}}
-\lambda_{9} \widehat{u}\, \overline{\widehat{\phi}}\right).
\end{eqnarray}
By multiplying \eqref{eq31}-\eqref{eq312} by $\lambda_1 ,\,-\lambda_2 ,\,\lambda_3 ,\,-\lambda_4$, $\lambda_5 ,\,1$ and 
$\lambda_6 ,\,\cdots,\,\lambda_{9}$, respectively, and adding the resulting equations, we find
\begin{eqnarray}\label{g26+23}
\frac{d}{dt}{F}_0 (\xi,\,t) &=& - \xi^{2}\left(k_3\lambda_3 \,\vert\widehat{\phi}\vert^{2} +\lambda_2\,\vert\widehat{u}\vert^{2} +(\lambda_5 -\lambda_1)\,\vert\widehat{y}\vert^{2} +(k_1\lambda_4 -k_1\lambda_5 -k_1 \lambda_2 )\,\vert\widehat{v}\vert^{2}\right)  \nonumber \\ 
&& - \xi^{2}\left(k_2\lambda_1\,\vert\widehat{z}\vert^{2} +k_4\,\vert\widehat{\sigma}\vert^{2} \right)
+I_1 Re (i\xi\widehat{\theta}\, \overline{\widehat{u}})+I_2 \xi^2 Re (\widehat{y}\, \overline{\widehat{\theta}}) +
\gamma\xi^{2} Re (\widehat{\sigma}\, \overline{\widehat{\theta}}) \nonumber \\ 
&& +\xi^{2} ((\lambda_3 +\lambda_4 )\vert\widehat{\theta}\vert^{2}
+\vert\widehat{\eta}\vert^{2}) +Re\left( \gamma\lambda_3\xi^2\,\widehat{\eta}\, \overline{\widehat{\phi}} - \gamma\lambda_{6}\xi^2\widehat{\eta}\, \overline{\widehat{z}} -ik_5 \xi^{2\epsilon_0 +1}\,\widehat{\eta}\, \overline{\widehat{\sigma}} -i\gamma\lambda_4\xi^3 \widehat{\eta}\, \overline{\widehat{v}}\right),
\end{eqnarray}
where 
\begin{equation*}
I_1 =\lambda_5\xi^2 +\lambda_2 -\lambda_{8}\quad\hbox{and}\quad I_2 = \lambda_4 -\lambda_5 - \lambda_6 - \lambda_7 .
\end{equation*}
We put  
\begin{equation*}
I_3 =\left(\frac{\vert\gamma\vert}{\gamma}k_4\lambda_0 +\gamma\right)\xi^2 +\frac{k_4}{\gamma}I_1\quad\hbox{and}\quad I_4 = \frac{k_4}{\gamma}(I_2\xi^2 +I_1 ), 
\end{equation*}
and introduce the functional
\begin{eqnarray}\label{F3}
F_1 (\xi,\,t) &=& \xi^2 F_0 (\xi,\,t)+\frac{\vert\gamma\vert}{\gamma}\lambda_0 \xi^2\,Re (i\xi \widehat{\theta}\, \overline{\widehat{\eta}})+\frac{k_4}{\gamma} I_1 \xi^2 Re\left(\widehat{v}\, \overline{\widehat{\sigma}}\right) 
+I_3 Re\left(i\xi\widehat{z}\, \overline{\widehat{\sigma}}\right)\nonumber \\ 
&& +\frac{1}{\gamma} I_1 \xi^2\,Re (\widehat{u}\, \overline{\widehat{\eta}})+I_4 Re\left(i\xi\widehat{\phi}\, \overline{\widehat{\sigma}}\right)-\frac{1}{\gamma} I_2 \xi^2 Re\left(i\xi\widehat{\eta}\, \overline{\widehat{y}}\right).
\end{eqnarray}
Multiplying \eqref{equp123}, \eqref{equp22}, \eqref{equp32}, \eqref{equp423}, \eqref{equp52}, \eqref{equp623} and \eqref{g26+23} by 
$\lambda_0 \xi^2$, $\frac{k_4}{\gamma} I_1$, $I_4$, 
$\frac{1}{\gamma} I_1 \xi^2$, $I_3$, $-\frac{1}{\gamma} I_2 \xi^2$ and $\xi^2$, respectively, and adding the obtained expressions, we arrive at (observe that \eqref{equp22}, \eqref{equp32} and \eqref{equp52} are valid also in case $(\tau_1 ,\tau_2 ,\tau_3)=(0,0,1)$)
\begin{eqnarray}\label{Fprime13}
\frac{d}{dt}{F}_1 (\xi,\,t) &=& - \xi^{4}\left(k_2\lambda_1\,\vert\widehat{z}\vert^{2} 
+k_3\lambda_3 \,\vert\widehat{\phi}\vert^{2} +\lambda_2\,\vert\widehat{u}\vert^{2} +\left(\lambda_5 -\lambda_1 \right)\,\vert\widehat{y}\vert^{2} +(k_1\lambda_5 -k_1\lambda_4 -k_1 \lambda_2 )\,\vert\widehat{v}\vert^{2}\right)  \nonumber \\ 
&& - \xi^{4}\left((\vert\gamma\vert\lambda_0 -\lambda_3 -\lambda_4)\vert\widehat{\theta}\vert^{2}+k_4\,\vert\widehat{\sigma}\vert^{2} \right)+(\vert\gamma\vert\lambda_0 +1 )\xi^{4} \vert\widehat{\eta}\vert^{2} \nonumber \\
&& +\xi^2 Re\left[\left(iI_5 \overline{\widehat{v}}+I_6 \overline{\widehat{\phi}}+I_7 \overline{\widehat{z}}-ik_5 \xi^{2\epsilon_0 +1}\overline{\widehat{\sigma}}+i\frac{\vert\gamma\vert}{\gamma}k_5 \lambda_0 \xi^{2\epsilon_0 +1}\overline{\widehat{\theta}}
+i\frac{k_5}{\gamma} \xi^{2\epsilon_0 +1}I_2 \overline{\widehat{y}} -\frac{k_5}{\gamma} \xi^{2\epsilon_0 }I_1 \overline{\widehat{u}}\right) \widehat{\eta}\right],
\end{eqnarray}
where
\begin{equation*}
I_5 =-\gamma\lambda_4\xi^3 +\left(\frac{\vert\gamma\vert}{\gamma}k_1\lambda_0 +\frac{k_4 -k_1}{\gamma} I_1 +\frac{k_1}{\gamma} I_2\right)\xi, 
\end{equation*}
\begin{equation*}
I_6 = \left(-\frac{\vert\gamma\vert}{\gamma} k_3 \lambda_0 +\gamma \lambda_3\right)\xi^2 + I_3 \quad\hbox{and}\quad 
I_7 =-\left(\frac{k_2}{\gamma}I_2 +\gamma\lambda_6 \right)\xi^2 +I_4 . 
\end{equation*}
We see that \eqref{I12} and \eqref{I57} are still valid. Then, applying Young's inequality, we get \eqref{I17}. Therefore, we define $F$ and $L$ by \eqref{FF1} and choose $0<\lambda_{3}$, $0<\lambda_{1} <\lambda_{4} <\lambda_{5}$,  $0<\lambda_{2} <\lambda_{5} -\lambda_{4}$, 
$\lambda_{0} >\frac{1}{\vert\gamma\vert}(\lambda_{3} +\lambda_{4})$ and 
\begin{equation*}
0<\varepsilon <\min\left\{k_2 \lambda_1 , k_3 \lambda_3 ,\lambda_2 ,\lambda_5 -\lambda_1 ,
k_1\lambda_5 -k_1\lambda_4 -k_1 \lambda_2 ,\vert\gamma\vert\lambda_0 -\lambda_3 -\lambda_4 , k_4 \right\}, 
\end{equation*}
we obtain \eqref{4g27+} and \eqref{4g30}. Consequentely, the proof can be ended as for Lemma \ref{lemma432}.
\end{proof}
\vskip0,1truecm
\begin{theorem}\label{theorem541}
The stability result given in Theorem \ref{theorem441} is satisfied when $(\tau_1 ,\tau_2 ,\tau_3)=(0,0,1)$.
\end{theorem}
\vskip0,1truecm
\begin{proof} 
The proof is identical to the one of Theorem \ref{theorem441}.
\end{proof}

\section{Application: lower order coupling terms \eqref{loct}}

This section concerns the stability of \eqref{e11} in case where the coupling terms \eqref{hoct} are replaced by the ones \eqref{loct}; more precisely, we study the stability of
\begin{align}\label{4s21}
\begin{cases}
\varphi_{tt} -k_1\,(\varphi_{x} +\psi +w)_{x} +\tau_1 \gamma q_{t}  = 0, \\
\psi_{tt} - k_2\,\psi_{xx} + k_1 \,(\varphi_{x} +\psi +w) + \tau_2 \gamma q_{t} = 0, \\
w_{tt} - k_3\,w_{xx} + k_1 (\varphi_{x} +\psi +w) +\tau_3 \gamma q_{t}  = 0,\\
q_{tt} -k_4 q_{xx} -k_5 q_{xxt} -\gamma (\tau_1 \varphi_{t} +\tau_2 \psi_{t} +\tau_3 w_{t})=0
\end{cases}
\end{align}
and 
\begin{align}\label{4s22}
\begin{cases}
\varphi_{tt} -k_1\,(\varphi_{x} +\psi +w)_{x} +\tau_1 \gamma q_{t}  = 0, \\
\psi_{tt} - k_2\,\psi_{xx} + k_1 \,(\varphi_{x} +\psi +w) + \tau_2 \gamma q_{t} = 0, \\
w_{tt} - k_3\,w_{xx} + k_1 (\varphi_{x} +\psi +w) +\tau_3 \gamma q_{t} = 0,\\
q_{tt} -k_4 q_{xx} +k_5 q_{t} -\gamma (\tau_1 \varphi_{t} +\tau_2 \psi_{t} +\tau_3 w_{t})=0
\end{cases}
\end{align}
with the initial conditions \eqref{i21}. We define $U$, its initial data $U_0$ and the energy $\widehat{E}$ as in Section 2. It is clear that \eqref{e11}, \eqref{g6}, \eqref{e14} and \eqref{Ep21} are valid with $A_2$ as in \eqref{e12},
\begin{equation}\label{4e12}
{A}_1 U_{x} = \left(
\begin{array}{c}
-u_{x} \\
\\
-k_1\,v_{x}  \\
\\
-\,y_{x}   \\
\\
-\,k_2 \,z_{x} \\
\\
-\,\theta_{x} \\
\\
-\,k_3\,\phi_{x} \\
\\
-\eta_x \\
\\
-k_4\sigma_x
\end{array}
\right) \,\,\hbox{and}\,\, A_0 U= \left(
\begin{array}{c}
-y - \theta \\
\\
\tau_1 \gamma\eta \\
\\
0   \\
\\
k_1\,v +\tau_2 \gamma\eta \\
\\
0  \\
\\
k_1\,v +\tau_3 \gamma\eta\\
\\
0 \\
\\
(1-\epsilon_0 )k_5\eta -\gamma (\tau_1 u +\tau_2 y +\tau_3 \theta ) 
\end{array}
\right)
\end{equation} 
So, instead of \eqref{e101}, we have
\begin{align}\label{4e101}
\begin{cases}
\widehat{v}_{t} - i\xi\widehat{u} - \widehat{y} - \widehat{\theta} = 0, \\
\widehat{u}_{t} - ik_1 \xi \,\widehat{v} +\tau_1 \gamma\,\widehat{\eta} = 0, \\
\widehat{z}_{t} - i\xi\widehat{y} = 0, \\
\widehat{y}_{t} - ik_2\xi\,\widehat{z} + k_1 \widehat{v} +\tau_2 \gamma\,\widehat{\eta} = 0, \\
\widehat{\phi}_{t} - i\xi\widehat{\theta} = 0,  \\
\widehat{\theta}_{t} - ik_3\,\xi\widehat{\phi} + k_1\,\widehat{v}+\tau_3 \gamma\,\widehat{\eta} = 0,\\
\widehat{\sigma}_{t} - i\xi\widehat{\eta} = 0,  \\
\widehat{\eta}_{t} -ik_4 \xi \widehat{\sigma}+k_5\xi^{2\epsilon_0}\widehat{\eta} -\gamma\, (\tau_1 \widehat{u} +\tau_2 \widehat{y}+\tau_3 \widehat{\theta}).
\end{cases}
\end{align}
\vskip0,1truecm
\begin{lemma}\label{lemma632}
Let $\widehat{U}$ be a solution of \eqref{g6}. Then, there exist $c,\,\widetilde{c}>0$ such that \eqref{11} holds true with the following $f$:
\vskip0,1truecm
{\bf Case $(\tau_1 ,\tau_2 ,\tau_3)=(1,0,0)$ under the condition $\chi \ne 0$}: 
\begin{align}\label{6ftildef}
f(\xi) =\frac{\xi^{4+2\epsilon_0 }}{{\tilde f}(\xi)}\quad\hbox{and}\quad {\tilde f}(\xi) = \begin{cases}
1 +\xi^{2} +\xi^4 +\xi^6 +\xi^{8} +\xi^{10} \quad &\hbox{for \eqref{4s21}},\\
1 +\xi^{2} +\xi^4 +\xi^6 +\xi^{8} \quad &\hbox{for \eqref{4s22}.}
\end{cases}
\end{align}
\vskip0,1truecm
{\bf Case $(\tau_1 ,\tau_2 ,\tau_3)=(1,0,0)$ under the condition $\chi = 0$}: $\vert\widehat{U}(\xi,\,t)\vert$ doesn't converge to zero when time $t$ goes to infinity.
\vskip0,1truecm
{\bf Case $(\tau_1 ,\tau_2 ,\tau_3)\in\{(0,1,0), (0,0,1)\}$}:
\begin{align}\label{7ftildef}
f(\xi) =\frac{\xi^{2 +2\epsilon_0 }}{{\tilde f}(\xi)}\quad\hbox{and}\quad {\tilde f}(\xi) = \begin{cases}
1 +\xi^{2} +\xi^{4} +\xi^{6} \quad &\hbox{for \eqref{4s21} under \eqref{k123}},\\
1 +\xi^{2} +\xi^{4} \quad &\hbox{for \eqref{4s22} under \eqref{k123}},\\
1 +\xi^{2} +\xi^4 +\xi^6 +\xi^{8} +\xi^{10} \quad &\hbox{for \eqref{4s21} without \eqref{k123}},\\
1 +\xi^{2} +\xi^4 +\xi^6 +\xi^{8} \quad &\hbox{for \eqref{4s22} without \eqref{k123}.}
\end{cases}
\end{align}
\end{lemma}
\vskip0,1truecm
\begin{proof} 
The proof is very similar to the one given in Sections 3 and 4 with some small modifications related to the coupling terms \eqref{loct}. We give here a bref idea of the proof.
\vskip0,1truecm
We see that, for \eqref{4e101}, the expressions \eqref{eq31}-\eqref{eq35} and \eqref{eq36}-\eqref{eq312} are satisfied with $\tau_j \gamma\widehat{\eta}$ instead of 
$i\tau_j \gamma\xi\widehat{\eta}$, and \eqref{eq31p} holds true if we replace 
$i\gamma\xi (\tau_1\widehat{u} +\tau_2\widehat{y}+\tau_3\widehat{\theta})$
by $-\gamma(\tau_1\widehat{u} +\tau_2\widehat{y}+\tau_3\widehat{\theta})$. 
\vskip0,1truecm
Now, we treat only the two cases  
$(\tau_1 ,\tau_2 ,\tau_3)=(1,0,0)$ and $(\tau_1 ,\tau_2 ,\tau_3)=(0,1,0)$. The last case $(\tau_1 ,\tau_2 ,\tau_3)=(0,0,1)$ can be treated using very similar modifications to the ones considered for the case $(\tau_1 ,\tau_2 ,\tau_3)=(0,1,0)$.

\subsection{Case $(\tau_1 ,\tau_2 ,\tau_3)=(1,0,0)$ under the condition $\chi \ne 0$}

We start by modifying the expressions \eqref{equ1}-\eqref{equ3p} (according to \eqref{4e101}). Multiplying \eqref{4e101}$_{2}$ and \eqref{4e101}$_{8}$ by $-\frac{\vert\gamma\vert}{\gamma}\,\xi^2\,\overline{\widehat{\eta}}$ and $-\frac{\vert\gamma\vert}{\gamma}\,\xi^2\,\overline{\widehat{u}}$, respectively, adding the resulting equations, taking the real part and using \eqref{hd1}, we get
\begin{eqnarray}\label{4equ1}
\frac{d}{dt}Re\left(-\frac{\vert\gamma\vert}{\gamma}\,\xi^2\,\widehat{u}\, \overline{\widehat{\eta}}\right) &=& \vert\gamma\vert\xi^{2}\left(\vert\widehat{\eta}\vert^{2} - \vert\widehat{u}\vert^{2}\right) - \frac{\vert\gamma\vert}{\gamma}k_4\xi^2\,Re\left(i\xi\widehat{\sigma}\,\overline{\widehat{u}}\right) \nonumber \\ 
&& -\frac{\vert\gamma\vert}{\gamma}k_1 \xi^2\,Re\left(i\xi\widehat{v}\, \overline{\widehat{\eta}}\right) +\frac{\vert\gamma\vert}{\gamma}k_5 \xi^{2\epsilon_0 +2}\,Re\left(\widehat{\eta}\, \overline{\widehat{u}}\right).
\end{eqnarray}
Multiplying \eqref{4e101}$_{6}$ and \eqref{4e101}$_{8}$ by $-i\xi\overline{\widehat{\eta}}$ and $i\xi\overline{\widehat{\theta}}$, respectively, adding the resulting equations, taking the real part and using \eqref{hd1}, we find
\begin{equation}\label{4equ2}
\frac{d}{dt}Re\left(i\xi\widehat{\eta}\, \overline{\widehat{\theta}}\right) = \gamma Re\left(i\xi\widehat{u}\,\overline{\widehat{\theta}}\right)-k_4\,\xi^2 Re\left(\widehat{\sigma}\, \overline{\widehat{\theta}}\right) -k_5 \xi^{2\epsilon_0}\,Re\left(i\xi\widehat{\eta}\, \overline{\widehat{\theta}}\right) +k_3 \xi^2\,Re\left(\widehat{\phi}\, \overline{\widehat{\eta}}\right)+k_1 Re\left(i\xi\widehat{v}\, \overline{\widehat{\eta}}\right).
\end{equation}
Also, multiplying \eqref{4e101}$_{4}$ and \eqref{4e101}$_{8}$ by $-i\xi\overline{\widehat{\eta}}$ and $i\xi\overline{\widehat{y}}$, respectively, adding the resulting equations, taking the real part and using \eqref{hd1}, we obtain
\begin{equation}\label{4equ3}
\frac{d}{dt}Re\left(i\xi\widehat{\eta}\, \overline{\widehat{y}}\right) = \gamma Re\left(i\xi\widehat{u}\,\overline{\widehat{y}}\right)-k_4 \xi^2\,Re\left(\widehat{\sigma}\, \overline{\widehat{y}}\right) -k_5 \xi^{2\epsilon_0}\,Re\left(i\xi\widehat{\eta}\, \overline{\widehat{y}}\right) +k_2 \xi^2\,Re\left(\widehat{z}\, \overline{\widehat{\eta}}\right)+k_1 Re\left(i\xi\widehat{v}\, \overline{\widehat{\eta}} \right).
\end{equation}
Multiplying \eqref{4e101}$_{1}$ and \eqref{4e101}$_{7}$ by $\overline{\widehat{\sigma}}$ and $\overline{\widehat{v}}$, respectively, adding the resulting equations, taking the real part and using \eqref{hd1}, we infer that
\begin{equation}\label{4equ1p}
\frac{d}{dt}Re\left(\widehat{v}\, \overline{\widehat{\sigma}}\right) = - Re\left(i\xi\widehat{\sigma}\,\overline{\widehat{u}}\right)+ Re\left(i\xi\widehat{v}\, \overline{\widehat{\eta}}\right) + Re\left(\widehat{y}\, \overline{\widehat{\sigma}}\right) + Re\left(\widehat{\theta}\, \overline{\widehat{\sigma}}\right).
\end{equation}
Similarily, multiplying \eqref{4e101}$_{3}$ and \eqref{4e101}$_{7}$ by $i\xi\overline{\widehat{\sigma}}$ and $-i\xi\overline{\widehat{z}}$, respectively, adding the resulting equations, taking the real part and using \eqref{hd1}, we arrive at
\begin{equation}\label{4equ2p}
\frac{d}{dt}Re\left(i\xi\widehat{z}\, \overline{\widehat{\sigma}}\right) = -\xi^2 Re\left( \widehat{\sigma}\,\overline{\widehat{y}}\right)+\xi^2 Re\left(\widehat{z}\, \overline{\widehat{\eta}}\right).
\end{equation}
Multiplying \eqref{4e101}$_{5}$ and \eqref{4e101}$_{7}$ by $i\xi\overline{\widehat{\sigma}}$ and $-i\xi\overline{\widehat{\phi}}$, respectively, adding the resulting equations, taking the real part and using \eqref{hd1}, we entail
\begin{equation}\label{4equ3p}
\frac{d}{dt}Re\left(i\xi\widehat{\phi}\, \overline{\widehat{\sigma}}\right) = -\xi^2 Re\left( \widehat{\sigma}\,\overline{\widehat{\theta}}\right)+\xi^2 Re\left(\widehat{\phi}\, \overline{\widehat{\eta}}\right).
\end{equation}
We put ${\tilde F}_0 (\xi,\,t)=\xi^2 {F}_0 (\xi,\,t)$, where ${F}_0$ is defined in \eqref{F0}. Multiplying \eqref{eq31}-\eqref{eq312} (with the modifications cited above) by 
$\lambda_1 ,\,\cdots ,\,\lambda_5$, $1$ and 
$\lambda_6 ,\,\cdots,\,\lambda_{9}$, respectively, and adding the obtained expressions, we find (instead of \eqref{g26+})
\begin{eqnarray}\label{4g26+}
\frac{d}{dt}{\tilde F}_0 (\xi,\,t) &=& - \xi^{4}\left(k_3\lambda_3 \,\vert\widehat{\phi}\vert^{2} +\left(\lambda_5 -\lambda_1 \right)\,\vert\widehat{y}\vert^{2} +\left(\lambda_4 -\lambda_3 \right)\,\vert\widehat{\theta}\vert^{2} +(k_1\lambda_2 -k_1\lambda_4 -k_1 \lambda_5 )\,\vert\widehat{v}\vert^{2}\right)  \nonumber \\ 
&& - \xi^{4}\left(k_2\lambda_1\,\vert\widehat{z}\vert^{2} 
+k_4\,\vert\widehat{\sigma}\vert^{2} \right)+I_1 \xi^{2} Re (i\xi\widehat{\theta}\, \overline{\widehat{u}})+I_2 \xi^{2} Re (i\xi\widehat{y}\, \overline{\widehat{u}}) -\gamma \xi^{2} Re (i\xi\widehat{\sigma}\, \overline{\widehat{u}}) \nonumber \\ 
&& +\xi^{4} (\lambda_2\vert\widehat{u}\vert^{2}
+\vert\widehat{\eta}\vert^{2}) +\xi^{2} Re\left( \gamma\lambda_8\widehat{\eta}\, \overline{\widehat{z}} + \gamma\lambda_{9}\widehat{\eta}\, \overline{\widehat{\phi}} -ik_5 \xi^{2\epsilon_0 +1}\,\widehat{\eta}\, \overline{\widehat{\sigma}} -i\gamma\lambda_2\xi \widehat{\eta}\, \overline{\widehat{v}}\right),
\end{eqnarray}
where $I_1$ and $I_2$ are defined in \eqref{46I1I2}. We put  
\begin{equation*}
I_3 = -\frac{k_4}{\gamma} (I_1 +\vert\gamma\vert\lambda_0 )\xi^2 -\gamma ,\quad I_4 = -\frac{k_4}{\gamma} (I_2 +\vert\gamma\vert\lambda_0 )\xi^2 -\gamma \quad\hbox{and}\quad I_5 = -\gamma -\frac{\vert\gamma\vert}{\gamma}k_4\lambda_0 \xi^2, 
\end{equation*}
and introduce the functional
\begin{eqnarray}\label{4F}
F_1 (\xi,\,t) &=& {\tilde F}_0 (\xi,\,t)-\frac{\vert\gamma\vert}{\gamma}\lambda_0 \xi^4\,Re (\widehat{u}\, \overline{\widehat{\eta}})+\frac{1}{\gamma} \xi^2 Re\left(iI_1\,\xi\widehat{\eta}\, \overline{\widehat{\theta}} +iI_2\xi\,\widehat{\eta}\, \overline{\widehat{y}}\right) \nonumber \\ 
&& +I_5 \xi^2\,Re (\widehat{v}\, \overline{\widehat{\sigma}})+I_4 Re\left(i\xi \widehat{z}\, \overline{\widehat{\sigma}}\right)+I_3 Re\left(i\xi \widehat{\phi}\, \overline{\widehat{\sigma}}\right).
\end{eqnarray}
Multiplying \eqref{4equ1}-\eqref{4equ3p} by $\lambda_0 \xi^2$, $\frac{1}{\gamma} I_1 \xi^2$, $\frac{1}{\gamma} I_2 \xi^2$, 
$I_5 \xi^2$, $I_4$ and $I_3$, respectively, adding the obtained equations, and adding \eqref{4g26+}, we arrive at
\begin{eqnarray}\label{4Fprime1}
\frac{d}{dt}{F}_1 (\xi,\,t) &=& - \xi^{4}\left(k_2\lambda_1\,\vert\widehat{z}\vert^{2} 
+k_3\lambda_3 \,\vert\widehat{\phi}\vert^{2} +\left(\lambda_5 -\lambda_1 \right)\,\vert\widehat{y}\vert^{2} +\left(\lambda_4 -\lambda_3 \right)\,\vert\widehat{\theta}\vert^{2} +(k_1\lambda_2 -k_1\lambda_4 -k_1 \lambda_5 )\,\vert\widehat{v}\vert^{2}\right)  \nonumber \\ 
&& - \xi^{4}\left((\vert\gamma\vert\lambda_0 -\lambda_2)\vert\widehat{u}\vert^{2}+k_4\,\vert\widehat{\sigma}\vert^{2} \right)+(\vert\gamma\vert\lambda_0 +1 )\xi^{4} \vert\widehat{\eta}\vert^{2} 
+\xi^2 Re\left( \frac{\vert\gamma\vert}{\gamma}k_5\lambda_0\xi^{2\epsilon_0 +2} \widehat{\eta}\, \overline{\widehat{u}} -ik_5 \xi^{2\epsilon_0 +1}\,\widehat{\eta}\, \overline{\widehat{\sigma}}\right) \nonumber \\
&& +\xi^2 Re\left[i\left(\frac{k_1}{\gamma}(I_1 +I_2 )+\gamma\lambda_2 +I_5 -\frac{\vert\gamma\vert}{\gamma} k_1\lambda_0 \xi^2\right)\xi\widehat{v}\, \overline{\widehat{\eta}} +i\frac{k_5}{\gamma}I_1 \xi^{2\epsilon_0 +1}\widehat{\theta}\, \overline{\widehat{\eta}}+i\frac{k_5}{\gamma}I_2 \xi^{2\epsilon_0 +1}\widehat{y}\, \overline{\widehat{\eta}}\right] \nonumber \\
&& +\xi^2 Re\left[\left(\gamma\lambda_8 +I_4 +\frac{k_2}{\gamma}I_2 \xi^2\right)\widehat{\eta}\, \overline{\widehat{z}} + \left(\gamma\lambda_{9} +I_3 +\frac{k_3}{\gamma}I_1 \xi^2\right)\widehat{\eta}\, \overline{\widehat{\phi}}\right].
\end{eqnarray}
Now, we consider ${\tilde f}$ and $f$ defined in \eqref{6ftildef}, and introduce the functionals  
\begin{eqnarray}\label{42FF1}
F (\xi,\,t) =\xi^{2\epsilon_0} F_1 (\xi,\,t)\quad\hbox{and}\quad L(\xi,\,t) =\lambda\,\widehat{E}(\xi,\,t) +
\frac{1}{{\tilde f}(\xi)}\,{F}(\xi,\,t).
\end{eqnarray}
Applying Young's inequality, \eqref{4Fprime1} implies \eqref{Fprime2}. So, the proof can be ended as for Lemma \ref{lemma32}. 

\subsection{Case $(\tau_1 ,\tau_2 ,\tau_3)=(1,0,0)$ under the condition $\chi = 0$}

To prove that $\vert\widehat{U}(\xi,\,t)\vert$ doesn't converge to zero when time $t$ goes to infinity, it is enough to prove
\eqref{AlambdaI}, where (according to \eqref{4e101}) 
\begin{equation*}
\lambda I -A = 
\begin{pmatrix}
\lambda & -i\xi & 0 & -1 & 0 & -1 & 0 & 0 \\
-ik_1 \xi & \lambda & 0 & 0 & 0 & 0 & 0 & \gamma \\
0 & 0 & \lambda & -i\xi & 0 & 0 & 0 & 0 \\
k_1 & 0 & -ik_2\xi & \lambda & 0 & 0 & 0 & 0 \\
0 & 0 & 0 & 0 & \lambda & -i\xi & 0 & 0 \\
k_1 & 0 & 0 & 0 & -ik_2 \xi & \lambda & 0 & 0 \\
0 & 0 & 0 & 0 & 0 & 0 & \lambda & -i\xi \\
0 & -\gamma & 0 & 0 & 0 & 0 & -ik_4\xi & k_5\xi^{2\epsilon_0} +\lambda
\end{pmatrix}.
\end{equation*}
A direct computaion shows that 
\begin{eqnarray*}
det(\lambda I -A) &=& 2k_1\lambda^2 (\lambda^2 +k_2 \xi^2 )\left[\lambda (\lambda +k_5 \xi^{2\epsilon_0})+k_4 \xi^2 +\gamma^2 \right] +k_4 \xi^2 (\lambda^2 +k_1 \xi^2 )\left(\lambda^2 + k_2\xi^2 \right)^2\\
&& +\lambda (\lambda^2 +k_2 \xi^2 )^2 \left[\lambda^2 (\lambda +k_5 \xi^{2\epsilon_0})+\gamma^2 \lambda +k_1 \xi^2 (\lambda +k_5 \xi^{2\epsilon_0})\right] . 
\end{eqnarray*}
Then, the same conclusions indicated in the proof of Theorem \ref{theorem410} are valide for \eqref{4e101}. 

\subsection{Case $(\tau_1 ,\tau_2 ,\tau_3)=(0,1,0)$}

First, we modify the expressions \eqref{equp12}-\eqref{equp62} according to \eqref{4e101}. 
Multiplying \eqref{4e101}$_{4}$ and \eqref{4e101}$_{8}$ by $-\frac{\vert\gamma\vert}{\gamma}\,\xi^2\,\overline{\widehat{\eta}}$ and $-\frac{\vert\gamma\vert}{\gamma}\,\xi^2\,\overline{\widehat{y}}$, respectively, adding the resulting equations, taking the real part and using \eqref{hd1}, we get
\begin{eqnarray}\label{4equp12}
\frac{d}{dt}Re\left(-\frac{\vert\gamma\vert}{\gamma}\,\xi^2\,\widehat{y}\, \overline{\widehat{\eta}}\right) &=& \vert\gamma\vert\xi^{2}\left(\vert\widehat{\eta}\vert^{2} - \vert\widehat{y}\vert^{2}\right) - \frac{\vert\gamma\vert}{\gamma}k_4\xi^2\,Re\left(i\xi\widehat{\sigma}\,\overline{\widehat{y}}\right) 
 +\frac{\vert\gamma\vert}{\gamma}k_1 \xi^2 \,Re\left(\widehat{\eta}\, \overline{\widehat{v}}\right) \nonumber \\ 
&& -\frac{\vert\gamma\vert}{\gamma}k_2 \xi^2\,Re\left(i\xi\widehat{z}\, \overline{\widehat{\eta}}\right) +\frac{\vert\gamma\vert}{\gamma}k_5 \xi^{2\epsilon_0 +2}\,Re\left(\widehat{\eta}\, \overline{\widehat{y}}\right).
\end{eqnarray}
Multiplying \eqref{4e101}$_{1}$ and \eqref{4e101}$_{7}$ by $i\xi\overline{\widehat{\sigma}}$ and $-i\xi\overline{\widehat{v}}$, respectively, adding the resulting equations, taking the real part and using \eqref{hd1}, we find
\begin{equation}\label{4equp22}
\frac{d}{dt}Re\left(i\xi\widehat{v}\, \overline{\widehat{\sigma}}\right) = - Re\left(i\xi\widehat{\sigma}\,\overline{\widehat{y}}\right)-\,Re\left(i\xi\widehat{\sigma}\, \overline{\widehat{\theta}}\right) -\xi^2 \,Re\left(\widehat{u}\, \overline{\widehat{\sigma}}\right) +\xi^2\,Re\left(\widehat{\eta}\, \overline{\widehat{v}}\right).
\end{equation}
Also, multiplying \eqref{4e101}$_{3}$ and \eqref{4e101}$_{7}$ by $\overline{\widehat{\sigma}}$ and $\overline{\widehat{z}}$, respectively, adding the resulting equations, taking the real part and using \eqref{hd1}, we obtain
\begin{equation}\label{4equp32}
\frac{d}{dt}Re\left(\widehat{z}\, \overline{\widehat{\sigma}}\right) = - Re\left(i\xi\widehat{\sigma}\,\overline{\widehat{y}}\right)+Re\left(i\xi\widehat{\eta}\, \overline{\widehat{z}}\right).
\end{equation}
Multiplying \eqref{4e101}$_{2}$ and \eqref{4e101}$_{8}$ by $-i\xi\overline{\widehat{\eta}}$ and $i\xi\overline{\widehat{u}}$, respectively, adding the resulting equations, taking the real part and using \eqref{hd1}, we infer that
\begin{equation}\label{4equp42}
\frac{d}{dt}Re\left(-i\xi\widehat{u}\, \overline{\widehat{\eta}}\right) = \gamma Re\left(i\xi\widehat{y}\,\overline{\widehat{u}}\right)+k_1 \xi^2 Re\left(\widehat{v}\, \overline{\widehat{\eta}}\right) -k_4 \xi^2 Re\left(\widehat{\sigma}\, \overline{\widehat{u}}\right) -k_5 \xi^{2\epsilon_0 }Re\left(i\xi\widehat{\eta}\, \overline{\widehat{u}}\right).
\end{equation}
Multiplying \eqref{4e101}$_{5}$ and \eqref{4e101}$_{7}$ by $\overline{\widehat{\sigma}}$ and $\overline{\widehat{\phi}}$, respectively, adding the resulting equations, taking the real part and using \eqref{hd1}, we infer that
\begin{equation}\label{4equp52}
\frac{d}{dt}Re\left(\widehat{\phi}\, \overline{\widehat{\sigma}}\right) = - Re\left(i\xi\widehat{\sigma}\,\overline{\widehat{\theta}}\right)+Re\left(i\xi\widehat{\eta}\, \overline{\widehat{\phi}}\right) .
\end{equation}
Finally, multiplying \eqref{4e101}$_{6}$ and \eqref{4e101}$_{8}$ by $\overline{\widehat{\eta}}$ and $\overline{\widehat{\theta}}$, respectively, adding the resulting equations, taking the real part and using \eqref{hd1}, it follows that
\begin{equation}\label{4equp62}
\frac{d}{dt}Re\left(\widehat{\eta}\, \overline{\widehat{\theta}}\right) = \gamma Re\left(\widehat{y}\,\overline{\widehat{\theta}}\right)+k_4 Re\left(i\xi\widehat{\sigma}\, \overline{\widehat{\theta}}\right) -k_3 Re\left(i\xi\widehat{\eta}\, \overline{\widehat{\phi}}\right) -k_5 \xi^{2\epsilon_0 }Re\left(\widehat{\eta}\, \overline{\widehat{\theta}}\right)-k_1 Re\left(\widehat{v}\, \overline{\widehat{\eta}}\right) .
\end{equation}
We define the functional ${F}_0$ by \eqref{F02}, and
we get (instead of \eqref{g26+2}) 
\begin{eqnarray}\label{4g26+2}
\frac{d}{dt}{F}_0 (\xi,\,t) &=& - \xi^{2}\left(k_3\lambda_3 \,\vert\widehat{\phi}\vert^{2} +\lambda_2\,\vert\widehat{u}\vert^{2} +\left(\lambda_4 -\lambda_3 \right)\,\vert\widehat{\theta}\vert^{2} +(k_1\lambda_5 -k_1\lambda_4 -k_1 \lambda_2 )\,\vert\widehat{v}\vert^{2}\right)  \nonumber \\ 
&& - \xi^{2}\left(k_2\lambda_1\,\vert\widehat{z}\vert^{2} 
+k_4\,\vert\widehat{\sigma}\vert^{2} \right)
+I_1 Re (i\xi\widehat{y}\, \overline{\widehat{u}})+I_2 \xi^2 Re (\widehat{y}\, \overline{\widehat{\theta}}) -
\gamma Re (i\xi\widehat{\sigma}\, \overline{\widehat{y}}) \nonumber \\ 
&& +\xi^{2} \left((\lambda_1 +\lambda_5)\vert\widehat{y}\vert^{2}
+\vert\widehat{\eta}\vert^{2}\right) +Re\left(-i \gamma\lambda_1\,\xi\widehat{\eta}\, \overline{\widehat{z}} +i \gamma\lambda_{7}\xi\widehat{\eta}\, \overline{\widehat{\phi}} -ik_5 \xi^{2\epsilon_0 +1}\,\widehat{\eta}\, \overline{\widehat{\sigma}} -\gamma\lambda_5\xi^2 \widehat{\eta}\, \overline{\widehat{v}}\right).
\end{eqnarray}
We put  
\begin{equation*}
I_3 =-\frac{\vert\gamma\vert}{\gamma}k_4\lambda_0 \xi^2 -\frac{k_4}{\gamma}I_1 -\gamma\quad\hbox{and}\quad I_4 = -\frac{k_4}{\gamma}(I_2\xi^2 +I_1 ), 
\end{equation*}
and introduce the functional
\begin{eqnarray}\label{4F}
F_1 (\xi,\,t) &=& F_0 (\xi,\,t)-\frac{\vert\gamma\vert}{\gamma}\lambda_0 \xi^2\,Re (\widehat{y}\, \overline{\widehat{\eta}})+\frac{k_4}{\gamma} I_1 Re\left(i\xi\widehat{v}\, \overline{\widehat{\sigma}}\right) 
+I_3 Re\left(\widehat{z}\, \overline{\widehat{\sigma}}\right)\nonumber \\ 
&& +\frac{1}{\gamma} I_1\,Re (i\xi\widehat{u}\, \overline{\widehat{\eta}})+I_4 Re\left(\widehat{\phi}\, \overline{\widehat{\sigma}}\right)-\frac{1}{\gamma} I_2 \xi^2 Re\left(\widehat{\eta}\, \overline{\widehat{\theta}}\right).
\end{eqnarray}
Multiplying \eqref{4equp12}-\eqref{4equp62} and \eqref{4g26+2} by $\lambda_0$, $\frac{k_4}{\gamma} I_1$, $I_3$, 
$-\frac{1}{\gamma} I_1$, $I_4$, $-\frac{1}{\gamma} I_2 \xi^2$ and $1$, respectively, and adding the obtained expressions, we arrive at 
\begin{eqnarray}\label{4FFprime1}
\frac{d}{dt}{F}_1 (\xi,\,t) &=& - \xi^{2}\left(k_2\lambda_1\,\vert\widehat{z}\vert^{2} 
+k_3\lambda_3 \,\vert\widehat{\phi}\vert^{2} +\lambda_2\,\vert\widehat{u}\vert^{2} +\left(\lambda_4 -\lambda_3 \right)\,\vert\widehat{\theta}\vert^{2} +(k_1\lambda_5 -k_1\lambda_4 -k_1 \lambda_2 )\,\vert\widehat{v}\vert^{2}\right)  \nonumber \\ 
&& - \xi^{2}\left((\vert\gamma\vert\lambda_0 -\lambda_1 -\lambda_5)\vert\widehat{y}\vert^{2}+k_4\,\vert\widehat{\sigma}\vert^{2} \right)+(\vert\gamma\vert\lambda_0 +1 )\xi^{2} \vert\widehat{\eta}\vert^{2} \nonumber \\
&& +\xi Re\left[\left(I_5 \overline{\widehat{v}}+iI_6 \overline{\widehat{z}}+iI_7 \overline{\widehat{\phi}}-ik_5 \xi^{2\epsilon_0}\overline{\widehat{\sigma}}+\frac{\vert\gamma\vert}{\gamma}k_5 \lambda_0 \xi^{2\epsilon_0 +1}\overline{\widehat{y}}
+\frac{k_5}{\gamma} \xi^{2\epsilon_0 +1}I_2 \overline{\widehat{\theta}} +i\frac{k_5}{\gamma} \xi^{2\epsilon_0}I_1 \overline{\widehat{u}}\right) \widehat{\eta}\right],
\end{eqnarray}
where
\begin{equation*}
I_5 =\left(\frac{\vert\gamma\vert}{\gamma}k_1\lambda_0 -\gamma\lambda_5 +\frac{k_1}{\gamma} (I_2 -I_1 ) +\frac{k_4}{\gamma} I_1 \right)\xi , 
\end{equation*}
\begin{equation*}
I_6 =\frac{\vert\gamma\vert}{\gamma}k_2\lambda_0 \xi^2 -\gamma\lambda_1 +I_3 \quad\hbox{and}\quad I_7 = \frac{k_3}{\gamma} I_2 \xi^2 + I_4 +\gamma \lambda_7 . 
\end{equation*}
Because \eqref{I12} is still satisfied, we infer that, for ${\tilde f}$ defined in \eqref{7ftildef},
\begin{eqnarray} \label{4FFprime2}
\frac{d}{dt}{F}_1 (\xi,\,t) &\leq& - \xi^{2}\left((k_2\lambda_1 -\varepsilon)\,\vert\widehat{z}\vert^{2} 
+(k_3\lambda_3 -\varepsilon)\,\vert\widehat{\phi}\vert^{2} +(\lambda_2 -\varepsilon)\,\vert\widehat{u}\vert^{2} +\left(\lambda_4 -\lambda_3 -\varepsilon\right)\,\vert\widehat{\theta}\vert^{2} \right)  \nonumber \\ 
&& - \xi^{2}\left((k_1\lambda_5 -k_1\lambda_4 -k_1 \lambda_2 -\varepsilon )\,\vert\widehat{v}\vert^{2} +(\vert\gamma\vert\lambda_0 -\lambda_1 -\lambda_5 -\varepsilon ) \vert\widehat{y}\vert^{2} +(k_4 -\varepsilon)\vert\widehat{\sigma} \vert^{2}\right) \nonumber \\
&& +C_{\varepsilon,\lambda_0 ,\cdots ,\lambda_9 } {\tilde f}(\xi)\vert\widehat{\eta}\vert^{2}.
\end{eqnarray}
Therefore, we introduce the functionals $F$ and $L$ defined in \eqref{FF1} and consider the same choices of $\lambda_0 ,\,\cdots ,\,\lambda_5$ and $\varepsilon$, we arrive at  
\begin{equation}\label{44g27+} 
\frac{d}{dt}{F}(\xi,\,t) \leq -c_1 \xi^{2+2\epsilon_0}\widehat{E} (\xi,\,t)+ C{\tilde f}(\xi)\xi^{2\epsilon_0} \vert\widehat{\eta}\vert^{2} .
\end{equation}
Hence, the proof can be finished as for Lemma \ref{lemma432}.
\end{proof}
\vskip0,1truecm
\begin{theorem}\label{theorem641}
Let $N,\,\ell\in \mathbb{N}^*$ such that $\ell\leq N$, $U_{0}\in H^{N}(\mathbb{R})\cap L^{1}(\mathbb{R})$
and $U$ be the solution of \eqref{e11}. Then, for any $j\in\{0,\,\ldots,\,N-\ell\}$, there exist $c_0 ,\,{\tilde c}_0 >0$ such that, for any $t\in \mathbb{R}_+$, 
\vskip0,1truecm
{\bf Case $(\tau_1 ,\tau_2 ,\tau_3)=(1,0,0)$ under the condition $\chi \ne 0$}:
\begin{equation*}
\Vert\partial_{x}^{j}U\Vert_{L^{2}(\mathbb{R})} \leq c_0 \,(1 + t)^{-1/12 - j/6}\,\Vert U_{0} \Vert_{L^{1}(\mathbb{R})} + c_0 \,(1 + t)^{-\ell/4}\,\Vert\partial_{x}^{j+\ell}U_{0} \Vert_{L^{2}(\mathbb{R})}
\end{equation*} 
for \eqref{4s21}, and
\begin{equation*}
\Vert\partial_{x}^{j}U\Vert_{L^{2}(\mathbb{R})} \leq c_0 \,(1 + t)^{-1/8 - j/4}\,\Vert U_{0} \Vert_{L^{1}(\mathbb{R})} + c_0 \,(1 + t)^{-\ell/4}\,\Vert\partial_{x}^{j+\ell}U_{0} \Vert_{L^{2}(\mathbb{R})} 
\end{equation*} 
for \eqref{4s22}.
\vskip0,1truecm
{\bf Case $(\tau_1 ,\tau_2 ,\tau_3)\in\{(0,1,0), (0,0,1)\}$ under the condition $k_1 =k_2 =k_3$}: 
\begin{equation*}
\Vert\partial_{x}^{j}U\Vert_{L^{2}(\mathbb{R})} \leq c_0 \,(1 + t)^{-1/8 - j/4}\,\Vert U_{0} \Vert_{L^{1}(\mathbb{R})} + c_0 (1 + t)^{-\ell/2}\,\Vert\partial_{x}^{j+\ell}U_{0} \Vert_{L^{2}(\mathbb{R})}
\end{equation*} 
for \eqref{4s21}, and
\begin{equation*}
\Vert\partial_{x}^{j}U\Vert_{L^{2}(\mathbb{R})} \leq c_0 \,(1 + t)^{-1/4 - j/2}\,\Vert U_{0} \Vert_{L^{1}(\mathbb{R})} + c_0 (1 + t)^{-\ell/2}\,\Vert\partial_{x}^{j+\ell}U_{0} \Vert_{L^{2}(\mathbb{R})}  
\end{equation*} 
for \eqref{4s22}.
\vskip0,1truecm
{\bf Case $(\tau_1 ,\tau_2 ,\tau_3)\in\{(0,1,0), (0,0,1)\}$ without the condition $k_1 =k_2 =k_3$}:
\begin{equation*}
\Vert\partial_{x}^{j}U\Vert_{L^{2}(\mathbb{R})} \leq c_0 \,(1 + t)^{-1/8 - j/4}\,\Vert U_{0} \Vert_{L^{1}(\mathbb{R})} + c_0 \,(1 + t)^{-\ell/6}\,\Vert\partial_{x}^{j+\ell}U_{0} \Vert_{L^{2}(\mathbb{R})}
\end{equation*}
for \eqref{4s21}, and
\begin{equation*}
\Vert\partial_{x}^{j}U\Vert_{L^{2}(\mathbb{R})} \leq c_0 \,(1 + t)^{-1/4 - j/2}\,\Vert U_{0} \Vert_{L^{1}(\mathbb{R})} + c_0 \,(1 + t)^{-\ell/6}\,\Vert\partial_{x}^{j+\ell}U_{0} \Vert_{L^{2}(\mathbb{R})} 
\end{equation*} 
for \eqref{4s22}.
\end{theorem}
\vskip0,1truecm
\begin{proof} 
The proof is identical to the one of Theorem \ref{theorem441}.
\end{proof}


\begin{thebibliography}{1}

\bibitem{agg} M. S. Alves, P. Gamboa, G. C. Gorain, A. Rambaud and O. Vera, Asymptotic behavior of a flexible structure with Cattaneo type of thermal effect, Indagationes Mathematicae, 27 (2016), 821-834.

\bibitem{apal1} T. A. Apalara, Uniform stability of a laminated beam with structural damping and second sound, ZAMP, 68 (2017), 40-55. 

\bibitem{apal2} T. A. Apalara, On the stability of thermoelastic laminated beam, Acta. Math. Scie., 39 (2019), 1517-1524. 

\bibitem{beim} C. F. Beards and I. M. A. Imam, The damping of plate vibration by interfacial slip between layers, Int. J. Mach. Tool. Des. Res., 18 (1978), 131-137. 

\bibitem{clx} X. G. Cao, D. Y. Liu and G. Q. Xu, Easy test for stability of laminated beams with structural
damping and boundary feedback controls, J. Dynamical Control Syst., 13 (2007), 313-336.

\bibitem{cdflr} M. M. Cavalcanti, V. N. Domingos Cavalcanti, F. A. Falcao Nascimento, I. Lasiecka and J. H. Rodrigues, Uniform decay rates for the energy of Timoshenko system with the arbitrary speeds of propagation and localized nonlinear damping, ZAMP, 65 (2014), 1189-1206.

\bibitem{chen} Z. Chen, W. Liu and D. Chen, General decay rates for a laminated beam with memory, Taiw. J. Math., 23 (2019), 1227-1252.

\bibitem{s6} L. Djouamai and B. Said-Houari, A new stability number of the Bresse-Cattaneo system, 
Math. Meth. Appl. Sci., 41 (2018), 2827-2847.

\bibitem{4} L. H. Fatori, R. N. Monteiro and H. D. Fern\a'{a}ndez Sare, The Timoshenko system with history and Cattaneo
law, Applied Mathematics and Computation, 228 (2014), 128-140.

\bibitem{feng} B. Feng, T. E. Ma, R. N. Monteiro and C. A. Raposo, Dynamics of laminated Timoshenko beams, J. Dyn. Diff. Equa., 30 (2018), 1489-1507.

\bibitem{s1} T. E. Ghoul, M. Khenissi and B. Said-Houari, On the stability of the Bresse system with frictional damping, J. Math. Anal. Appl., 455 (2017), 1870-1898.

\bibitem{gues2} A. Guesmia, Asymptotic stability of Bresse system with one infinite memory in the longitudinal displacements, Medi. J. Math., 14 (2017), 19 pages.

\bibitem{gues4} A. Guesmia, Non-exponential and polynomial stability results of a Bresse system with one infinite memory in the vertical displacement, Nonauton. Dyn. Syst., 4 (2017), 78-97.

\bibitem{gues8} A. Guesmia, Well-posedness and stability results for laminated Timoshenko beams with interfacial slip and infinite memory, IMA J. Math. Cont. Info., 37 (2020), 300-350.

\bibitem{gms} A. Guesmia, S. Messaoudi and A. Soufyane, On the stabilization for a linear Timoshenko system with infinite history and applications to the coupled Timoshenko-heat systems, Elec. J. Diff. Equa., 2012 (2012), 1-45.

\bibitem{hans} S. W. Hansen, In control and estimation of distributed parameter systems: Non-linear
phenomena, International Series of Numerical Analysis, 118 (1994), 143-170.

\bibitem{hasp} S. W. Hansen and R. Spies, Structural damping in a laminated beams due to interfacial slip,
J. Sound Vibration, 204 (1997), 183-202.

\bibitem{8} K. Ide, K. Haramoto and S. Kawashima, Decay property of regularity-loss type for dissipative Timoshenko
system, Math. Mod. Meth. Appl. Sci., 18 (2008), 647-667.

\bibitem{s4} M. Khader and B. Said-Houari, Decay rate of solutions to Timoshenko system with past history in unbounded domains, 
Appl. Math. Optim., 75 (2017), 403-428.

\bibitem{s5} M. Khader and B. Said-Houari, Optimal decay rate of solutions to Timoshenko system with past history in unbounded domains, Z. Anal. Anwend, 37 (2018), 435-459.

\bibitem{li} G. Li, X. Kong and W. Liu, General decay for a laminated beam with structural damping and memory: the case of non-equal wave speeds, J. Inte. Equa., 30 (2018), 95-116.

\bibitem{lizh} W. Liu and W. Zhao, Exponential and polynomial decay for a laminated beam with Fourier's type heat conduction, Preprints 2017, 2017020058, doi: 10.20944/preprints201702.0058.v1.

\bibitem{lota1} A. Lo and N. E Tatar, Stabilization of laminated beams with interfacial slip, Elec. J. Diff. Equa., 2015 (2015), 1-14.

\bibitem{lota2} A. Lo and N. E. Tatar, Uniform stability of a laminated beam with structural memory, Qual. Theory Dyn.
Syst., 15 (2016), 517-540.

\bibitem{lota3} A. Lo and N. E. Tatar, Exponential stabilization of a structure with interfacial slip, Discrete Contin. Dyn.
Syst., 36 (2016), 6285-6306.

\bibitem{must} M. I. Mustafa, Laminated Timoshenko beams with viscoelastic damping, J. Math. Anal. Appl., 466 (2018), 619-641.

\bibitem{rapo} C. A. Raposo, Exponential stability for a structure with interfacial slip and
frictional damping, Appl. Math. Lett., 53 (2016), 85-91.

\bibitem{rvma} C. A. Raposo, O. V. Villagr\'an, J. E. Mu\~noz Rivera and M. S. Alves, Hybrid laminated Timoshenko beam, J. Math. Phys., 58 (2017), 11 pages. 

\bibitem{s3} B. Said-Houari and R. Racke, Decay rates and global existence for semilinear dissipative Timoshenko systems, 
Quart. Appl. Math., 71 (2013), 229-266.

\bibitem{6} B. Said-Houari and R. Rahali, Asymptotic behavior of the Cauchy problem of the Timoshenko system
in thermoelsaticity of type III, Evol. Equa. Cont. Theory, 2 (2013), 423-440.

\bibitem{s2} B. Said-Houari and A. Soufyane, The effect of frictional damping terms on the decay rate of the Bresse system, Evol. Equa. Cont. Theory, 3 (2014), 713-738.

\bibitem{10} M. L. Santos, D. S. Almeida and J. E. Mu\~noz Rivera, The stability number of the Timoshenko system
with second sound, J. Diff. Equa., 253 (2012), 2715-2733.

\bibitem{7} A. Soufyane and B. Said-Houari, The effect of the wave speeds and the frictional damping terms on the
decay rate of the Bresse system, Evol. Equa. Cont. Theory, 3 (2014), 713-738.

\bibitem{tesc} G. Teschl, Ordinary differential equations and dynamical systems, Amer. Math.
Soc., 140 (2012), ISBN 978-0-8218-8328-0.

\bibitem{tata} N. E. Tatar, Stabilization of a laminated beam with interfacial slip by boundary controls, Boundary Value Problem,
2015, DOI: 10.1186/s13661-015-0432-3.  

\bibitem{wxy}  J. M. Wang, G. Q. Xu and S. P. Yung, Exponential stabilization of laminated beams with structural damping and boundary feedback controls, SIAM J. Control Optim., 44 (2005), 1575-1597.
\end{thebibliography}
\end{document}